\documentclass{amsart}
\usepackage{amsmath}
\usepackage{amsthm}
\usepackage{amsfonts}
\usepackage{amssymb}

\newlength{\abstand}
\def\Aut{{\text{Aut}}}
\def\Autk{{\Aut_k}}
\def\AutkL{{\Aut_kL}}
\def\AutN{{\Aut\,N}}

\def\F{{\mathbb F}}

\def\C{{\mathbb C}}
\def\N{{\mathbb N}}
\def\Q{{\mathbb Q}}
\def\Z{{\mathbb Z}}
\def\P{{\mathbb P}}
\def\Cat{{\mathcal C}}
\def\Ck{{{\mathcal C}_k}}
\def\CK{{{\mathcal C}_K}}
\def\Ckp{{{\mathcal C}'_k}}
\def\CKp{{{\mathcal C}'_K}}

\def\Fkz{{{\bf\mathcal F}_k^{n,m}}}

\def\FKz{{{\bf\mathcal F}_K^{n,m}}}
\def\FPk{{\widetilde{{\bf\mathcal F}_k^{n,r}}}}
\def\FPK{{\widetilde{{\bf\mathcal F}_K^{n,r}}}}

\def\FPz{{\widetilde{{\bf\mathcal F}_k^{n,m}}}}
\def\FPZ{{\widetilde{{\bf\mathcal F}_K^{n,m}}}}
\def\FX{{F_{\mathcal X}}}
\def\FXs{{F^*_{\mathcal X}}}
\def\GLnK{{\text{GL}(n,K)}}

\def\kranz{{S_n{\textstyle\int}\mu_m}}
\def\kranzN{{\mu_m^N\hspace{-1mm}\rtimes\hspace{-1mm}\AutN}}
\def\Kt{{K^\times}}
\def\KtcN{{K^\times_c}^N}
\def\Lt{{L^\times}}
\def\LtLtm{{L^\times/{L^\times}^m}}
\def\mum{{\mu_m}}
\def\mumn{{\mu_m^n}}
\def\mumN{{\mu_m^N}}
\def\Pmn{{P^m_n}}

\def\Vark{{{\bf Var}_k}}
\def\VarK{{{\bf Var}_K}}
\def\VarkI{{{\bf Var}_k^{\text{Iso}}}}
\def\VarKI{{{\bf Var}_K^{\text{Iso}}}}
\def\Ratk{{{\bf Var}'_k}}
\def\RatK{{{\bf Var}'_K}}

\def\MC{{\bf Rep}_\Cat^S}
\def\MCk{{\bf Rep}_\Cat^{G_k}}
\def\MCK{{\bf Rep}_\Cat^{G_K}}
\def\Mk{{\bf Rep}_L^{G_k}}
\def\MK{{\bf Rep}_L^{G_K}}

\def\PGLnK{{\text{PGL}(n,K)}}
\def\Qlk{{\bf Rep}_{\Q_l}^{G_k}}
\def\QlK{{\bf Rep}_{\Q_l}^{G_K}}
\def\Mot{{\boldsymbol{\mathcal M}}}
\def\MotEk{{\boldsymbol{{\mathcal M}^E_k}}}
\def\Hc{{H_{\text{cont}}^1}}
\def\Zc{{Z_{\text{cont}}^1}}
\def\Het{{\text{H}_{\text{\'{e}t}}^i}}
\def\Hetn{{\text{H}_{\text{\'{e}t}}^{n-2}}}
\def\Hetnull{{\text{H}_{\text{\'{e}t}}^0}(\bar{\mathcal X},\Q_l)}
\def\Heteins{{\text{H}_{\text{\'{e}t}}^1}(\bar{\mathcal X},\Q_l)}
\def\Hetzwei{{\text{H}_{\text{\'{e}t}}^2}(\bar{\mathcal X},\Q_l)}

\def\homkLK{{\text{Hom}_k(L,K)}}

\def\fata{\mbox{\boldmath $a$}}
\def\fatap{\mbox{\boldmath $a'$}}
\def\fatb{\mbox{\boldmath $b$}}
\def\fatc{\mbox{\boldmath $c$}}
\def\fatbp{\mbox{\boldmath $b'$}}
\def\fatcp{\mbox{\boldmath $c'$}}

\def\resp{resp.\ }
\def\sign{{\text{sgn}}}
\def\star{\hspace{-1mm}*\hspace{-1mm}}
\def\tA{{\tilde{A}}}
\def\zmzs{{(\Z/m\Z)^\times}}
\newcommand\op[2]{{{}^{#1}\hspace{-0.7 mm}{#2}}}
\newcommand\spec[1]{{\text{Spec}\,#1}}

\swapnumbers
\theoremstyle{definition}
\newtheorem{defi}{Definition}[section]
\newtheorem{bsp}[defi]{Example}
\newtheorem{bspe}[defi]{Examples}
\newtheorem{satzdefi}[defi]{Proposition/ Definition}
\newtheorem{lemma}[defi]{Lemma}
\newtheorem{lemmadefi}[defi]{Lemma/ Definition}
\newtheorem{bem}[defi]{Remark}
\newtheorem{satz}[defi]{Proposition}
\newtheorem{thm}[defi]{Theorem}
\newtheorem{cor}[defi]{Corollary}
\newtheorem{reminder}[defi]{Reminder}

\input xy
\xyoption{all}

\title{On the Zeta Function of Forms of Fermat Equations}
\author[L. Br\"unjes]{Lars Br\"unjes}
\date{\today}

\thanks{The author is supported by the European Union.}
\thanks{This work has been written under hospitality of the Department of Pure
  Mathematics and Mathematical Statistics of the University of Cambridge}
\address{University of Cambridge \\ Department of Pure Mathematics and
  Mathematical Statistics \\ Centre for Mathematical Sciences \\ Wilberforce
  Road \\ Cambridge CB3 0WB \\ United Kingdom}
\email{lbrunjes.gmx.de}
\date{\today}

\begin{document}

\begin{abstract}
  We study ``forms of the Fermat equation'' over an arbitrary field $k$,
  i.e. homogenous equations of degree $m$ in $n$ unknowns that can be transformed 
  into the Fermat equation $X_1^m+\ldots+X_n^m$ by a suitable linear change of variables
  over an algebraic closure of $k$.
  Using the method of Galois descent, we classify all such forms.
  In the case that $k$ is a finite field of characteristic greater than $m$ 
  that contains the $m$-th roots of unity,
  we compute the Galois representation on $l$-adic cohomology 
  (and so in particular the zeta function) of the hypersurface associated to an
  arbitrary form of the Fermat equation. 
\end{abstract}

\maketitle

%%%%%%%%%%%%%%%%%%%%%%%%%%%%%%%%%%%%%%%%%%%%%%%%%%%%%%%%%%%%%%%%%%%%%%%%%%%%%%%%%%%%%%%%%%%%%%%%%%%%%%%%%%%%%%
%%%%%%%%%%%%%%%%%%%%%%%%%%%%%%%%%%%%%%%%%%%%%%%%%%%%%%%%%%%%%%%%%%%%%%%%%%%%%%%%%%%%%%%%%%%%%%%%%%%%%%%%%%%%%%
%%%%%%%%%%%%%%%%%%%%%%%%%%%%%%%%%%%%%%%%%%%%%%%%%%%%%%%%%%%%%%%%%%%%%%%%%%%%%%%%%%%%%%%%%%%%%%%%%%%%%%%%%%%%%%
%%%%%%%%%%%%%%%%%%%%%%%%%%%%%%%%%%%%%%%%%%%%%%%%%%%%%%%%%%%%%%%%%%%%%%%%%%%%%%%%%%%%%%%%%%%%%%%%%%%%%%%%%%%%%%
%%%%%%%%%%%%%%%%%%%%%%%%%%%%%%%%%%%%%%%%%%%%%%%%%%%%%%%%%%%%%%%%%%%%%%%%%%%%%%%%%%%%%%%%%%%%%%%%%%%%%%%%%%%%%%
%%%%%%%%%%%%%%%%%%%%%%%%%%%%%%%%%%%%%%%%%%%%%%%%%%%%%%%%%%%%%%%%%%%%%%%%%%%%%%%%%%%%%%%%%%%%%%%%%%%%%%%%%%%%%%
%%%%%%%%%%%%%%%%%%%%%%%%%%%%%%%%%%%%%%%%%%%%%%%%%%%%%%%%%%%%%%%%%%%%%%%%%%%%%%%%%%%%%%%%%%%%%%%%%%%%%%%%%%%%%%
%%%%%%%%%%%%%%%%%%%%%%%%%%%%%%%%%%%%%%%%%%%%%%%%%%%%%%%%%%%%%%%%%%%%%%%%%%%%%%%%%%%%%%%%%%%%%%%%%%%%%%%%%%%%%%
%%%%%%%%%%%%%%%%%%%%%%%%%%%%%%%%%%%%%%%%%%%%%%%%%%%%%%%%%%%%%%%%%%%%%%%%%%%%%%%%%%%%%%%%%%%%%%%%%%%%%%%%%%%%%%

\section*{Introduction}

%\markboth{INTRODUCTION}{}

\vspace{\abstand}

One of the most thoroughly studied equations in algebraic geometry is the
\emph{Fermat equation}
\[
  \Pmn:X_1^m+X_2^m+\ldots+X_n^m=0
\]
for natural numbers $m,n\geq 2$ and unknowns $X_1,\ldots,X_n$ in a ring $R$.
In the case $m\geq 3$, $n=3$ and $R=\Z$,
it is the subject of \emph{Fermat's Last Theorem} which was proved by Andrew Wiles in 1993.
In the case $R=\F_q$ a finite field, it served Andr\'{e} Weil in his important paper \cite{weil}
as motivation to formulate his famous \emph{Weil conjectures}, whose proof was only completed in 1973 by
Pierre Deligne and without doubt marked one of the highlights of 20th century mathematics.

Also in present day algebraic geometry, $\Pmn$ (\resp the \emph{Fermat hypersurface} $X^m_n$ 
defined by $\Pmn$)
is often used as an example or test case,
e.g. in the study of the conjectures of Hodge and Tate.
And even though --- mainly because of Tetsuji Shioda's work
(\cite{shi1}, \cite{shioda79b}, \cite{shioda79}, \cite{shioda81}, \cite{shioda82}, \cite{shi2}, 
\cite{shi3}, \cite{shi4}) ---
a lot more is known about $X^m_n$ than about hypersurfaces in general,
a lot of questions remain open; for example, Shioda was only able to prove the Hodge and the Tate conjecture
for a lot of Fermat hypersurfaces, but not for all of them.\\

For a homogenous polynomial $P$ over $k:=\F_q$ in $n\geq 2$ unknowns with associated projective variety $X$,
one of the basic invariants one can look at is its \emph{zeta function} which is defined as\\[-4 mm]
\[
  \zeta(P,t):=\zeta(X,t):=\text{exp}\left(\sum_{i=1}^\infty\frac{\#X(\F_{q^i})}{i}t^i\right)
  \in \Q(t)\subset\Q((t))
\]
and which can be computed using

\begin{thm}\label{satzzeta}(Deligne \cite{sga7ii})\\
  If $X$ is \emph{regular}, we have
  \[
    \zeta(P,t)={\underbrace{Q(P,t)}_{\in\Z[t]}}
    ^{\left[(-1)^{n+1}\right]}\prod_{i\in\{0,\ldots,n-2\}\setminus\{\frac{n-2}{2}\}}
    \frac{1}{1-q^it},
  \]
  and if $l\nmid q$ is a prime,
  $\bar{X}:=X\times_{\F_q}\bar{\F}_q$
  and $F:\bar{X}\rightarrow\bar{X}$ is the geometric Frobenius, we have
  \[
    Q(P,t)=\text{det}\,\left(1-F^*t\,\vert\,\Hetn(\bar{X},\Q_l)\right).
  \]
\end{thm}

The zeta function of $X^m_n$
was already known to Weil,
and Deligne \cite{deligne} computed the Frobenius action on
$H^*_{\text{\'{e}t}}(\bar{X}^m_n,\Q_l)$, which in this case is fairly simple,
because the cohomology decomposes into canonical "motivic" 1-dimensional subspaces on whom
the Frobenius action is given by special Hecke characters, the so called \emph{Jacobi Sums}
(see theorem \ref{satzkohom} below).\\

Starting from Fermat hypersurfaces, it is a natural step to next consider the slightly more general
class of \emph{diagonal hypersurfaces}, given by \emph{diagonal equations} of the form
\[
  a_1X_1^m+\ldots+a_nX_n^m=0
\]
for constants $a_1,\ldots,a_n\in k$;
these are studied intensely (in the case $q\equiv 1\pmod{m}$) by
Fernando Q. Gouv\^{e}a und Noriko Yui in \cite{yui3},
which makes it for example easy to compute their zeta function.\\

The $\bar{k}$-linear substitution $X_i\mapsto\sqrt[m]{a_i}X_i$ transforms the diagonal equation into
the Fermat equation $\Pmn$,
which means that the two equations are isomorphic \emph{over $\bar{k}$} or, as we say,
that the diagonal equation is a \emph{form} of $\Pmn$.

This naturally leads to the questions whether there are other forms of $\Pmn$ that are not diagonal.\\

Using the method of Galois descent (\ref{satzdefitheta}), we will see that forms of $\Pmn$ are
classified by the (non-abelian) cohomology $\Hc(G_k,(\kranz)/\mum)$, and it turns out that the diagonal
equations are exactly those that already come from $\Hc(G_k,\mumn)$, which means that they form a rather
special subclass of all forms.

Our aim in this paper therefore is to consider \emph{all} forms of $\Pmn$ at the same time, using Galois descent to first
classify them and then to compute their zeta function.\\

Other than in the case of diagonal equations, it is in general difficult to see whether a given equation
is a form of $\Pmn$. An example is the following equation
\begin{equation}\label{eqbsp1}
  \begin{array}{lcrl}
    0 & = & & 2a^3+6a^2c+a^2d+6ab^2+2abc+4abd+2ac^2+3acd+6ad^2 \\
    & & & \mbox{}\hspace{15mm}+5b^3+2b^2c+5b^2d+5bc^2+5bcd+6bd^2 \\
    & & & \mbox{}\hspace{15mm}+2c^3+6c^2d+cd^2+d^3 \\
    & & + & y^3+4x^2y+5y^3
  \end{array}
\end{equation}
\mbox{}\\[-3 mm]
over $\F_7$ which is a form of $P^3_6$.\\

The plan of this paper is as follows:
In the first chapter we introduce the general notion of a \emph{coefficient extension}
which is an axiomatization of situations in which Galois descent works,
something we then explain in the second chapter.
In the third chapter we apply our machinery to the problem of classifying all forms of $\Pmn$
over an arbitrary field $k$,
and in the fourth chapter we compute the zeta function of all such forms in the case $k=\F_q$,
$k^\times\supseteq\mum$. We sum everything up in our main result --- Theorem \ref{thmzeta} ---
and conclude in the last chapter with some remarks and open questions on the case $k^\times\not\supseteq\mum$.

\vspace{\abstand}

%%%%%%%%%%%%%%%%%%%%%%%%%%%%%%%%%%%%%%%%%%%%%%%%%%%%%%%%%%%%%%%%%%%%%%%%%%%%%%%%%%%%%%%%%%%%%%%%%%%%%%%%%%%%%%
%%%%%%%%%%%%%%%%%%%%%%%%%%%%%%%%%%%%%%%%%%%%%%%%%%%%%%%%%%%%%%%%%%%%%%%%%%%%%%%%%%%%%%%%%%%%%%%%%%%%%%%%%%%%%%
%%%%%%%%%%%%%%%%%%%%%%%%%%%%%%%%%%%%%%%%%%%%%%%%%%%%%%%%%%%%%%%%%%%%%%%%%%%%%%%%%%%%%%%%%%%%%%%%%%%%%%%%%%%%%%
%%%%%%%%%%%%%%%%%%%%%%%%%%%%%%%%%%%%%%%%%%%%%%%%%%%%%%%%%%%%%%%%%%%%%%%%%%%%%%%%%%%%%%%%%%%%%%%%%%%%%%%%%%%%%%
%%%%%%%%%%%%%%%%%%%%%%%%%%%%%%%%%%%%%%%%%%%%%%%%%%%%%%%%%%%%%%%%%%%%%%%%%%%%%%%%%%%%%%%%%%%%%%%%%%%%%%%%%%%%%%
%%%%%%%%%%%%%%%%%%%%%%%%%%%%%%%%%%%%%%%%%%%%%%%%%%%%%%%%%%%%%%%%%%%%%%%%%%%%%%%%%%%%%%%%%%%%%%%%%%%%%%%%%%%%%%
%%%%%%%%%%%%%%%%%%%%%%%%%%%%%%%%%%%%%%%%%%%%%%%%%%%%%%%%%%%%%%%%%%%%%%%%%%%%%%%%%%%%%%%%%%%%%%%%%%%%%%%%%%%%%%
%%%%%%%%%%%%%%%%%%%%%%%%%%%%%%%%%%%%%%%%%%%%%%%%%%%%%%%%%%%%%%%%%%%%%%%%%%%%%%%%%%%%%%%%%%%%%%%%%%%%%%%%%%%%%%
%%%%%%%%%%%%%%%%%%%%%%%%%%%%%%%%%%%%%%%%%%%%%%%%%%%%%%%%%%%%%%%%%%%%%%%%%%%%%%%%%%%%%%%%%%%%%%%%%%%%%%%%%%%%%%

\section{Coefficient Extensions}

\vspace{\abstand}

Let $K/k$ be an arbitrary Galois extension of fields with Galois group
$G:=\mbox{Gal}\,(K/k)$.\\[\abstand]

\begin{defi}\label{deficoeffext}
A {\em coefficient extension} (from $k$ to $K$)
consists of two categories $\Ck$ and $\CK$,
a covariant functor $F:\Ck\rightarrow\CK$
and a $G$-action (from the left) on $\text{Iso}_{\,\CK}(FY,FZ)$ for all
$Y,Z\in\text{Ob}(\Ck)$, so that the following two conditions are satisfied:
\begin{enumerate}
  \item[(CE1)]
    The action is compatible with compositions, i.e.
    for objects $X,Y,Z\in\text{Ob}\,(\Ck)$,
    iso\-morph\-isms
    $X\xrightarrow{g}Y$ and $Y\xrightarrow{f}Z$
    and an element $s\in G$ we have:
    \[
      \op{s}{(fg)}=\op{s}{f}\op{s}{g}.
    \]
  \item[(CE2)]
    Exactly those isomorphisms that come from $\Ck$ are fix under the action of $G$,
    i.e. for objects $Y,Z\in\text{Ob}(\Ck)$ we have:
    \[
      \text{Im}\left(\text{Iso}_{\,\Ck}(Y,Z)\xrightarrow{F}\text{Iso}_{\,\CK}(FY,FZ)\right)
      \ =\
      \Bigl[\text{Iso}_{\,\CK}(FY,FZ)\Bigr]^G.
    \]
\end{enumerate}
\end{defi}

\vspace{\abstand}

\begin{lemma}\label{lemmakestrich}
  Let $\Ck$, $\CK$ and $F:\Ck\longrightarrow\CK$ be as in \ref{deficoeffext}, and for
  all objects $X$ and $Y$ from $\Ck$ let us have an action of $G$ on $\text{Mor}_K(FX,FY)$,
  such that the conditions (CE1') and (CE2') that we get by replacing ``isomorphism'' by
  ``morphism'' in (CE1) and (CE2) are satisfied.\\[1 mm]
  Then we get a coefficient extension after restricting the $G$-action to isomorphisms.
\end{lemma}

\vspace{\abstand}

\begin{proof}[Proof:]
  This is obvious, seeing as (CE1') and (CE2') imply immediately
  that the orbit of an isomorphism under the $G$-action
  only consists of isomorphisms.
\end{proof}

\vspace{\abstand}

\begin{bspe}\label{bspkoeff}
\mbox{}\\[-5 mm]
\begin{enumerate}
  \item\label{bspkoeffi}
    Let $\Vark$ and $\VarK$ be the categories of quasiprojective varieties over $k$ \resp $K$,
    let $F$ be the functor of base change with $K$ over $k$,
    and let the action of $s\in G$ be defined by the commutativity of the following diagram
    \begin{equation}\label{eq1}
      \xymatrix@R=10mm@C=40mm{
        Y_K \ar@{}[rd]|{=} \ar[r]^f & Z_K \\
        Y_K \ar[u]^{1_Y\times\spec{s}}_\wr \ar@{-->}[r]_{\op{s}{f}} &
          Z_K \ar[u]^\wr_{1_Z\times\spec{s}}
      }
    \end{equation}
    Then descent theory (\cite[5.3.]{sga1}) shows that $\Vark\xrightarrow{F}\VarK$
    is a coefficient extension.
  \item\label{bspkoeffrat}
    Let $\Ratk$ and $\RatK$ be the categories of \emph{geometrically irreducible}
    quasiprojective varieties over $k$ \resp $K$ with \emph{dominant, rational maps} as morphisms.
    Define $\Ratk\xrightarrow{F}\RatK$ as in \ref{bspkoeffi} as base change,
    and again define the $G$-action by the commutativity of \eqref{eq1}.
    Then $\Ratk\xrightarrow{F}\RatK$ is a coefficient extension.
  \item\label{bspkoeffiii}
    For $n,m\in\N_+$ consider the categories $\Fkz$ and $\FKz$
    whose objects are non-zero homogenous polynomials of degree $m$
    in $X_1,\ldots,X_n$ over $k$ \resp $K$ and
    whose morphisms $P\rightarrow Q$ are elements $A=(a_{ij})$
    of $\text{GL}(n,k)$ (\resp $\GLnK)$,
    with $Q(AX):=Q\left(\sum_{j=1}^na_{1j}X_j,\ldots,\sum_{j=1}^na_{nj}X_j\right)=P$.
    If for $F:\Fkz\rightarrow\FKz$ we take the obvious functor
    and for an (iso-)morphism $(a_{ij})\in\GLnK$ in $\FKz$ and $s\in G$ we put
    $\op{s}{(a_{ij})}:=(\op{s}{a_{ij}})$,
    we get a coefficient extension $\Fkz\xrightarrow{F}\FKz$.\\[4 mm]
    For example, if $\text{char}\,(k)\not\in\{2,3\}$, and we consider
    $P:={X_1^3+X_2^3+X_3^3}$ and
    $Q:={X_1^3+\frac{1}{12}X_2^2X_3+\frac{1}{4}X_3^3}$
    in ${\bf\mathcal F}_k^{3,3}$,
    then
    $A:=\Bigl(\begin{array}{rrr}
      \scriptstyle 0 & \scriptstyle  0 & \scriptstyle 1\\[-3 mm]
      \scriptstyle 3 & \scriptstyle -3 & \scriptstyle 0\\[-3 mm]
      \scriptstyle 1 & \scriptstyle  1 & \scriptstyle 0
    \end{array}\Bigr)$
    defines a morphism $P\rightarrow Q$ because of
    $Q(X_3,3X_1-3X_2,X_1+X_2)=P$.
  \item\label{bspkoeffiv}
    For $n,m\in\N_+$ now consider the categories $\FPz$ and $\FPZ$
    whose objects are those of $\Fkz$ \resp $\FKz$
    but whose morphisms $P\rightarrow Q$ are elements of $\text{PGL}(n,k)$ (\resp $\PGLnK)$,
    represented by regular $n\times n$-matrices $A=(a_{ij})$
    with $Q(AX)=\lambda\cdot P$
    for a $\lambda\in k^\times$ (\resp $\Kt$).
    For $F:\FPz\rightarrow\FPZ$ we again take the obvious functor,
    and for an (iso-)morphism $(\overline{a_{ij}})\in\PGLnK$ in $\FPZ$ and $s\in G$ we put
    $\op{s}{\ \overline{(a_{ij})}}:=\overline{(\op{s}{a_{ij}})}$
    and in this way get another coefficient extension $\FPz\xrightarrow{F}\FPZ$.\\[4 mm]
    In the explicit example from \ref{bspkoeffiii}, $\bar{A}\in\text{PGL}(3,k)$ now defines a
    morphism $P\rightarrow Q':=12Q=12X_1^3+X_2^2X_3+3X_3^3$
    in $\widetilde{{\bf\mathcal F}_k^{3,3}}$.
  \item\label{bspkoeffv}
    For a category $\Cat$ and a group $S$ let $\MC$ be the category
    whose objects are pairs $(X,\varphi)$
    with $X\in\text{Ob}(\Cat)$ and $\varphi$ an $S$-action on $X$,
    and whose morphisms are $S$-equivariant morphisms of $\Cat$.\\[2 mm]
    Let $G_k$ and $G_K$ be the absolute Galois groups of $k$ and $K$, so that we have
    $G=G_k/G_K$.
    Let $\Cat$ be an arbitrary category, and consider the categories
    $\MCk$ and $\MCK$ and the functor $F:\MCk\rightarrow\MCK, (X,\varphi)\mapsto(X,\varphi\vert G_K)$.
    For $s\in G$ choose a representant $\tilde{s}\in G_k$ and define the action of $s$
    on an isomorphism $(X,\varphi\vert G_K)\xrightarrow{f}(Y,\psi\vert G_K)$
    by $\op{s}{f}:=\psi(\tilde{s})\circ f\circ\varphi(\tilde{s})^{-1}$.
    Then $\MCk\xrightarrow{F}\MCK$ is a coefficient extension.\\[2 mm]
    (In the special case that $\Cat$ is the category of vector spaces over a field $L$,
    we write $\Mk$ and $\MK$ instead of $\MCk$ and $\MCK$.)
  \item\label{bspkoeffvi}
    For a category $\mathcal C$ let ${\mathcal C}^{\text{Iso}}$ denote the category
    whose objects are those of $\mathcal C$ and whose morphisms are the \emph{isomorphisms} of $\mathcal C$.
    If $\Ck\xrightarrow{F}\CK$ is a coefficient extension, then clearly
    ${\mathcal C}_k^{\text{Iso}}\xrightarrow{F\vert{\mathcal C}_k^{\text{Iso}}}{\mathcal C}_K^{\text{Iso}}$
    is one as well.
\end{enumerate}
\end{bspe}

\vspace{\abstand}

%%%%%%%%%%%%%%%%%%%%%%%%%%%%%%%%%%%%%%%%%%%%%%%%%%%%%%%%%%%%%%%%%%%%%%%%%%%%%%%%%%%%%%%%%%%%%%%%%%%%%%%%%%%%%%
%%%%%%%%%%%%%%%%%%%%%%%%%%%%%%%%%%%%%%%%%%%%%%%%%%%%%%%%%%%%%%%%%%%%%%%%%%%%%%%%%%%%%%%%%%%%%%%%%%%%%%%%%%%%%%
%%%%%%%%%%%%%%%%%%%%%%%%%%%%%%%%%%%%%%%%%%%%%%%%%%%%%%%%%%%%%%%%%%%%%%%%%%%%%%%%%%%%%%%%%%%%%%%%%%%%%%%%%%%%%%
%%%%%%%%%%%%%%%%%%%%%%%%%%%%%%%%%%%%%%%%%%%%%%%%%%%%%%%%%%%%%%%%%%%%%%%%%%%%%%%%%%%%%%%%%%%%%%%%%%%%%%%%%%%%%%
%%%%%%%%%%%%%%%%%%%%%%%%%%%%%%%%%%%%%%%%%%%%%%%%%%%%%%%%%%%%%%%%%%%%%%%%%%%%%%%%%%%%%%%%%%%%%%%%%%%%%%%%%%%%%%
%%%%%%%%%%%%%%%%%%%%%%%%%%%%%%%%%%%%%%%%%%%%%%%%%%%%%%%%%%%%%%%%%%%%%%%%%%%%%%%%%%%%%%%%%%%%%%%%%%%%%%%%%%%%%%
%%%%%%%%%%%%%%%%%%%%%%%%%%%%%%%%%%%%%%%%%%%%%%%%%%%%%%%%%%%%%%%%%%%%%%%%%%%%%%%%%%%%%%%%%%%%%%%%%%%%%%%%%%%%%%
%%%%%%%%%%%%%%%%%%%%%%%%%%%%%%%%%%%%%%%%%%%%%%%%%%%%%%%%%%%%%%%%%%%%%%%%%%%%%%%%%%%%%%%%%%%%%%%%%%%%%%%%%%%%%%
%%%%%%%%%%%%%%%%%%%%%%%%%%%%%%%%%%%%%%%%%%%%%%%%%%%%%%%%%%%%%%%%%%%%%%%%%%%%%%%%%%%%%%%%%%%%%%%%%%%%%%%%%%%%%%

\section{Forms}

\vspace{\abstand}

\begin{reminder}\label{definonabel}
  Let $S$ be a (topological) group and $A$ a (discrete) $S$-group, i.e. a (topological) group
  on which $S$ acts (continuously) via group-homomorphisms from the left.
  We briefly want to remind ourselves of some definitions concerning \emph{non-abelian cohomology}
  (for details refer to \cite{serre} or \cite{larsphd}):
  \begin{enumerate}
    \item
      The \emph{zeroth cohomology (of $S$ with values in $A$)} is the group
      \[
        H^0(S,A):=A^S:=\left\{a\in A\,\vert\,\forall s\in S\,:\,\op{s}a=a\right\}.
      \]
      $H^0(S,\_)$ obviously becomes a covariant functor
      from the category of (discrete) $S$-groups to the category of groups.
    \item
      $Z^1(S,A)$ (\resp $Z_{\text{cont}}^1(S,A)$),
      the set of \emph{(continuous) 1-cocycles from $S$ in $A$}, is the pointed set
      of (continuous) maps $s\mapsto a_s$ from $S$ to $A$,
      satisfying the usual 1-cocycle-condition $a_{st}=a_s\op{s}{a_t}$,
      where the special point is given by the \emph{trivial 1-cocycle} $s\mapsto 1$.
      On $Z^1(S,A)$ (\resp $\Zc(S,A)$),
      we have the equivalence relation of being \emph{cohomologous}, given by
      \[
        (a_s)\sim(a'_s)\ :\Leftrightarrow\ \exists b\in A:\forall s\in S: a'_s=b^{-1}a_s\op{s}{b},
      \]
      and the pointed set $H^1(S,A)$ (\resp $\Hc(S,A)$) of equivalence classes is called
      the {\em first (continuous) cohomology of $S$ with values in $A$}.
      It is easy to see that $H^1(S,\_)$ (\resp $\Hc(S,\_)$)
      induces a covariant functor from the category of (discrete) $S$-groups to the category of pointed sets.
    \item\label{defitwist}
      Let $1\rightarrow A\xrightarrow{i}B\xrightarrow{p}C\rightarrow 1$ be an exact sequence of (discrete)
      $S$-groups, and let $b=(b_s)$ be a (continuous) 1-cocycle from $S$ in $B$. Then we can define
      new (continuous) $S$-actions on $A$, $B$ and $C$ by \emph{twisting by $b$}, i.e. by setting
      \[
        \op{s}{\,'a}:=i^{-1}\left(b_s\cdot\op{s}{(ia)}\cdot b_s^{-1}\right),\;\;\;\;
        \op{s}{\,'b}:=b_s\cdot\op{s}{b}\cdot b_s^{-1},\;\;\;\;
        \op{s}{\,'c}:=(pb_s)\cdot\op{s}{c}\cdot(pb_s)^{-1},
      \]
      for $s\in S$, $a\in A$, $b\in B$ and $c\in C$. These $S$-groups are denoted by $A_b$, $B_b$ and $C_b$,
      and we get a bijection
      \[
        \begin{array}{ccc}
          H^1_{\text{(cont)}}(S,A_b)/H^0(S,C_b) &
          \stackrel{\sim}{\longrightarrow} &
          \bigl\{\bigl.[b']\in H^1_{\text{(cont)}}(S,B)\,\bigr\vert\,p_*[b']=p_*[b]\bigr\} \\
          (a_s) &
          \mapsto &
          (a_s\cdot b_s)
        \end{array}
      \]
      where the right-$H^0(S,C_b)$-action on $H^1_{\text{(cont)}}(S,A_b)$ is given as follows:
      \[
        \begin{array}{rclcl}
          H^1_{\text{(cont)}}(S,A_b) &
          \times &
          H^0(S,C_b) &
          \longrightarrow &
          H^1_{\text{(cont)}}(S,A_b) \\[2 mm]
          \bigl((a_s) &
          , &
          c\bigr) &
          \mapsto &
          (\beta^{-1}a_sb_s\op{s}{\beta}b_s^{-1})
        \end{array}
      \]
      \begin{flushright}
        (for a $\beta$ with $p(\beta)=c$).
      \end{flushright}
  \end{enumerate}
\end{reminder}

\vspace{\abstand}

\begin{satzdefi}\label{satzdefitheta}(Galois descent)
  Let $F:\Ck\rightarrow\CK$ be a coefficient extension and $X\in\text{Ob}(\Ck)$.
  We define the class of {\em $\CK/\Ck$-forms of $X$} as
  \[
    E(\CK/\Ck,X)\ :=\ \{Y\in\text{Ob}(\Ck)\,|\,FY\cong_KFX\}/\text{$\Ck$-isomorphisms}.
  \]
  According to (CE1), the group $A(X):=\Aut_K(FX)$ is a $G$-group, and
  if for $[Y]\in E(\CK/\Ck,X)$ we choose $f:FY\xrightarrow{\sim}FX$,
  then a straight forward computation shows that
  $\left[s\mapsto f\op{s}{(f^{-1})}\right]$ is a 1-cocycle of $G$ in $A(X)$
  whose class $\vartheta[Y]$ in $H^1(G,A(X))$ only depends on $[Y]$ and that
  the correspondence $[Y]\mapsto\vartheta[Y]$ is
  \emph{injective}.
  In particular, $E(\CK/\Ck,X)$ is a \emph{set},
  and we thus get a well-defined injective map (of pointed sets)
  \[
    \fbox{$\vartheta:\ E(\CK/\Ck,X)\hookrightarrow H^1(G,A(X))$}.
  \]
\end{satzdefi}

\vspace{\abstand}

\begin{satz}\label{corautform}
  Let $\Ck\xrightarrow{F}\CK$ be a coefficient extension, $X\in\text{Ob}(\Ck)$ and
  $Y\in E(\CK/\Ck,X)$ with $\vartheta[Y]=[a]$ for a 1-cocycle $a$ from $G$ in $A(X)$.
  If $\Aut_k(Y)\xrightarrow{F}A(Y)$ is
  injective\footnote[2]{This condition is automatically satisfied for all coefficient extensions
  in \ref{bspkoeff}.},
  then
  \[
    \fbox{$\Aut_k(Y)\cong H^0\left(G,A(X)_a\right)$.}
  \]
\end{satz}

\vspace{\abstand}

\begin{proof}
  Let $f:FY\xrightarrow{\sim}FX$ be an isomorphism with $a_s=f\op{s}{(f^{-1})}$. Then it is easy to see
  that $A(Y)\rightarrow A(X)_a$, $b\mapsto fbf^{-1}$ is an isomorphism of $G$-groups. Thus the proposition
  follows from $(CE2)$.
\end{proof}

\vspace{\abstand}

\begin{lemmadefi}\label{lemmadefistetig}
  Let $F:\Ck\rightarrow\CK$ be a coefficient extension.
  If for all $X,Y\in\text{Ob}\,(\Ck)$ the $G$-action on
  $\text{Iso}_K\,(FX,FY)$ is \emph{continuous}
  (where we consider $G$ and $\text{Iso}_K\,(FX,FY)$
  as topological spaces, endowed with Krull \resp discrete topology),
  we call that coefficient extension \emph{continuous}.\\[2 mm]
  In that case one sees easily that
  $\vartheta$ maps $E(\CK/\Ck,X)$ into $H^1_{\text{cont}}(G,A(X))$.
\end{lemmadefi}

\vspace{\abstand}

\begin{bsp}\label{bspstetig}
  It is easy to see that the coefficient extensions
  $\Vark\xrightarrow{F}\VarK$, $\Ratk\xrightarrow{F}\RatK$, $\Fkz\xrightarrow{F}\FKz$
  and $\FPz\xrightarrow{F}\FPZ$ from \ref{bspkoeff} are continuous,
  whereas the coefficient extension $\MCk\xrightarrow{F}\MCK$ in general is not.
\end{bsp}

\vspace{\abstand}

\begin{lemmadefi}\label{ldeftwist}
  Let $L$ be an arbitrary field,
  $F:\Mk\rightarrow\MK$ the coefficient extension from \ref{bspkoeff}\ref{bspkoeffv},
  $V:=(V,\varphi)\in\text{Ob}(\Mk)$ and
  $\xi=(a_{\bar{s}})\in Z^1(G,A(V))$.\\[2 mm]
  We define a new $G_k$-action $\varphi^\xi$ on $V$ by $\varphi^\xi(s):=a_{\bar{s}}\varphi(s)$.
  Then it is easy to see that $V(\xi):=(V,\varphi^\xi)$ is a $K/k$-form of $V$ with
  \[
    \vartheta[V(\xi)]\ \ =\ \ [\xi]\in H^1(G,A(V)).
  \]
  In particular, we see that $\vartheta$ is bijective in this case.
\end{lemmadefi}

\vspace{\abstand}

\begin{defi}\label{defiMKE}
  A \emph{ morphism of coefficient extensions} $\Ck\xrightarrow{F}\CK$ and
  $\Ckp\xrightarrow{F'}\CKp$
  is given by \emph{covariant} functors
  $H_k:\Ck\rightarrow\Ckp$ and $H_K:\CK\rightarrow\CKp$
  and an isomorphism of functors $h:H_KF\xrightarrow{\sim}F'H_k$,
  such that for all $Y,Z\in\text{Ob}(\Ck)$ the composition
  \[
    \text{Iso}_{\CK}(FY,FZ)
    \xrightarrow{H_K}
    \text{Iso}_{\CK'}(H_KFY,H_KFZ)
    \xrightarrow{h}
    \text{Iso}_{\CK'}(F'H_kY,F'H_kZ)
  \]
  is $G$-equivariant.
\end{defi}

\vspace{\abstand}

The following proposition is an easy consequence of our definitions:

\vspace{\abstand}

\begin{satz}\label{satzMKE}
  If $\Ck\xrightarrow{F}\CK\xrightarrow{(H_k,H_K,h)}\Ckp\xrightarrow{F'}\CKp$
  is a morphism of coefficient extensions,
  then for all $X\in\text{Ob}(\Ck)$ the following diagram of pointed sets commutes:\\[2 mm]
  \[\xymatrix@R=20mm@C=25mm{
    {E(\CK/\Ck,X)} \ar[rr]^{\textstyle\varphi}_{\textstyle[Y]\mapsto[H_kY]}
      \ar@{^{(}->}[d]_{\textstyle\vartheta} \ar@{}[rrd]|{=}
      & & {E(\CKp/\Ckp,H_kX)} \ar@{^{(}->}[d]^{\textstyle\vartheta} \\
    {H^1(G,A(X))} \ar[rr]_{\textstyle\psi}^{\textstyle\left(a_s\right)\mapsto\left(hH_Ka_s\right)}
      & & {H^1(G,A(H_kX))}
  }\]
\end{satz}

\vspace{\abstand}

\begin{bspe}\label{bspMKE}
\mbox{}\\[-5 mm]
\begin{enumerate}
  \item
     For natural numbers $n\geq 2$ and $r\geq 1$,
     there is an obvious morphism of coefficient extensions
     $\bigl(F:\FPk\rightarrow\FPK\bigr)\longrightarrow\bigl(F':\VarkI\rightarrow\VarKI\bigr)$.
  \item
    For $l\neq\text{char}(k)$ prime and $i\in\N_0$,
    $l$-adic cohomology induces a morphism of coefficient extensions
    $\bigl(F':\VarkI\rightarrow\VarKI\bigr)\longrightarrow\bigl(F'':\Qlk\rightarrow\QlK\bigr)$
    by sending $X\in\text{Ob}(\VarkI)$ to $\Het(X\times_k\bar{K},\Q_l)$
    and an isomorphism $f$ to $(f^{-1})^*$, analogous for $\VarKI$.
  \item
    Combining (i) and (ii), we get the following morphism of coefficient extensions:\\[2 mm]
    \hspace*{2 cm}
    \xymatrix@C=1.5 cm@R=3 mm{
      \FPK \ar[rrr]^{H_K} & & & \QlK \\
      & & & \\
      & \ar@{=>}[ru]^{h} & & \\
      \FPk \ar[rrr]_{H_k}\ar[uuu]^{F} & & & \Qlk \ar[uuu]_{F''}
    }
\end{enumerate}
\end{bspe}

\vspace{\abstand}

\begin{cor}\label{corladischeKohom}
Let $n\geq 2$, $r\geq 1$ and $i\geq 0$ be natural numbers,
$l\neq\text{char}(k)$ prime,
$P\in\text{Ob}(\FPk)$,
$X\subset\P_k^{n-1}$ the hypersurface defined by $P$
and $\varphi_P$ the action of $G_k$ on $l$-adic cohomology:
\[
  \varphi_P:G_k\longrightarrow\Aut_{\Q_l}\,\underbrace{\Het(X\times_k\bar{K},\Q_l)}_{=:V}.
\]
Take a $\FPK/\FPk$-Form $Q$ of $P$ with
$\vartheta[Q]=:[(A_s)]\in\text{H}^1(G,A(P))$ and corresponding hypersurface $Y$.
Then $\varphi_Q$, the corresponding $G_k$-action on $V$ (note $Y_K\cong X_K$!), is given by
\begin{equation}\label{eqopallg}
  \boxed{
  \varphi_Q:G_k\longrightarrow\Aut_{\Q_l}(V),\
  s\mapsto\Het\left(A_{\bar{s}}^{-1}\right)\circ\varphi_P(s)
  }
\end{equation}
In the particular case that $k$ is \emph{finite} with arithmetic Frobenius $f\in G_k$, we get:
\begin{equation}\label{eqopfin}
  \boxed{
    F_Y^*=\Het\left(\op{\bar{f}^{-1}}{A_{\bar{f}}}\right)\circ F_X^*
  }
\end{equation}
where $F_P^*:=\varphi_P(f)^{-1}$ \resp $F_Q^*:=\varphi_Q(f)^{-1}$ denotes geometric Frobenii,
and the knowledge of $F_Q^*$ determines $\varphi_Q$.
\end{cor}

\vspace{\abstand}

\begin{proof}
\eqref{eqopallg} and \eqref{eqopfin} are easy consequences of
\ref{ldeftwist} and \ref{satzMKE}, and the last statement follows from
\ref{lemmadefistetig} and \ref{bspstetig}
and the fact that $f$ generates $G_k$ topologically.
\end{proof}

\vspace{\abstand}

%%%%%%%%%%%%%%%%%%%%%%%%%%%%%%%%%%%%%%%%%%%%%%%%%%%%%%%%%%%%%%%%%%%%%%%%%%%%%%%%%%%%%%%%%%%%%%%%%%%%%%%%%%%%%%
%%%%%%%%%%%%%%%%%%%%%%%%%%%%%%%%%%%%%%%%%%%%%%%%%%%%%%%%%%%%%%%%%%%%%%%%%%%%%%%%%%%%%%%%%%%%%%%%%%%%%%%%%%%%%%
%%%%%%%%%%%%%%%%%%%%%%%%%%%%%%%%%%%%%%%%%%%%%%%%%%%%%%%%%%%%%%%%%%%%%%%%%%%%%%%%%%%%%%%%%%%%%%%%%%%%%%%%%%%%%%
%%%%%%%%%%%%%%%%%%%%%%%%%%%%%%%%%%%%%%%%%%%%%%%%%%%%%%%%%%%%%%%%%%%%%%%%%%%%%%%%%%%%%%%%%%%%%%%%%%%%%%%%%%%%%%
%%%%%%%%%%%%%%%%%%%%%%%%%%%%%%%%%%%%%%%%%%%%%%%%%%%%%%%%%%%%%%%%%%%%%%%%%%%%%%%%%%%%%%%%%%%%%%%%%%%%%%%%%%%%%%
%%%%%%%%%%%%%%%%%%%%%%%%%%%%%%%%%%%%%%%%%%%%%%%%%%%%%%%%%%%%%%%%%%%%%%%%%%%%%%%%%%%%%%%%%%%%%%%%%%%%%%%%%%%%%%
%%%%%%%%%%%%%%%%%%%%%%%%%%%%%%%%%%%%%%%%%%%%%%%%%%%%%%%%%%%%%%%%%%%%%%%%%%%%%%%%%%%%%%%%%%%%%%%%%%%%%%%%%%%%%%
%%%%%%%%%%%%%%%%%%%%%%%%%%%%%%%%%%%%%%%%%%%%%%%%%%%%%%%%%%%%%%%%%%%%%%%%%%%%%%%%%%%%%%%%%%%%%%%%%%%%%%%%%%%%%%
%%%%%%%%%%%%%%%%%%%%%%%%%%%%%%%%%%%%%%%%%%%%%%%%%%%%%%%%%%%%%%%%%%%%%%%%%%%%%%%%%%%%%%%%%%%%%%%%%%%%%%%%%%%%%%

\section{Forms of the Fermat equation}

\vspace{\abstand}

\label{notation}
Consider $m\in\N_{\geq 3}$, $n\in\N_+$,
$k$ a field with $(\text{char}\,k)\nmid m!$,
$K$ a separable algebraic closure of $k$,
$G:=G_k=\text{Gal}(K/k)$ the absolute Galois group of $k$ and
$\Pmn:=X_1^m+\ldots+X_n^m\in k[X_i]$ the \emph{Fermat equation}. The latter is an object of both
$\Fkz$ and $\FPz$.

\vspace{\abstand}

\begin{satz}(Shioda \cite{shi4})\label{satzshioda}\\
  Let $\kranz$ denote the Kranz product of the symmetric group $S_n$ and the group $\mum$ of
  $m$-th roots of unity in $K$, i.e. the semidirect product $\mumn\rtimes S_n$.
  Then
  \[
    \kranz\longrightarrow\GLnK,\;\;\;
    (\mu_1,\ldots,\mu_n)\cdot\sigma\mapsto(\mu_i\cdot\delta_{i,\sigma(j)})_{ij}
  \]
  defines an isomorphism from $\kranz$ onto $\Aut_\Fkz(\Pmn)$.
\end{satz}

\vspace{\abstand}

\begin{reminder}\label{remindersym}(see for example \cite{tamme})\\
  Let $L$ be a commutative, separable $k$-algebra of degree $n$. If $N$ denotes the set $\homkLK$ and
  if we choose a bijection $N\xrightarrow{\sim}\{1,\ldots,n\}$
  and thereby identify the groups $\AutN$ and $S_n$,
  then $s\mapsto[\varphi\mapsto s\circ\varphi]$ defines a continuous
  1-cocycle of $G$ in $N$ \resp $S_n$.
  If we start with an isomorphic $k$-algebra or choose another bijection, we end up with a
  cohomologous 1-cocycle and therefore get a well defined map from
  isomorphism classes of commutative, separable
  $k$-algebras of degree $n$ to $\Hc(G,S_n)$. This map is an isomorphism of pointed sets.
\end{reminder}

\vspace{\abstand}

\begin{thm}\label{thmkranz}(Rupprecht \cite{chris}, Br\"unjes \cite{larsphd})
  \begin{enumerate}
    \item\label{kranziso}
      The canonical projection $\kranz\stackrel{p}{\twoheadrightarrow}S_n$ induces a surjection
      $\Hc(G,\kranz)$ $\stackrel{p_*}\twoheadrightarrow\Hc(G,s_n)$ on cohomology.
      Let $c\in\Zc(G,S_n)$ be an arbitrary 1-cocycle,
      represented in the sense of \ref{remindersym} by a commutative, separable $k$-Algebra
      $L$ of degree $n$.
      If $N$ denotes $\homkLK$
      and if we choose a bijection $N\xrightarrow{\sim}\{1,\ldots,n\}$ and thereby identify
      $\kranz$ with $\kranzN$,
      we have a bijection
      \begin{equation}\label{eqcocycle}\fbox{$
        \begin{array}{ccc}
          \AutkL\backslash(\LtLtm) & \xrightarrow{\sim} &
          \left\{\left.[b]\in\Hc(G,\kranzN)\ \right\vert\ p_*[b]=[c]\right\} \\[2 mm]
          \bar{x} & \mapsto &
          \left(\left(\frac{s\sqrt[m]{s^{-1}\varphi x}}{\sqrt[m]{\varphi x}}\right)_{\varphi\in N}
          \cdot c_s\right)_s
        \end{array}
      $}\end{equation}
      and therefore an isomorphism
      \begin{equation}\label{eqisohc}
        \boxed{
          \coprod_{[L]}
          \AutkL\backslash \LtLtm
          \stackrel{\sim}{\longrightarrow}
          \Hc\bigl(G,A(\Pmn)\bigr)
        }
      \end{equation}
      of pointed sets where $[L]$ runs over the isomorphism classes of commutative, separable
      $k$-algebras of degree $n$.
    \item\label{kranztwist}
      The map $E(\FKz/\Fkz,\Pmn)\xrightarrow{\vartheta}\Hc\bigl(G,A(\Pmn)\bigr)$ is an
      \emph{isomorphism}. If $(L,x)$ represents an element of the right hand side via \eqref{eqisohc}
      and if $\{e_1,\ldots,e_n\}$ is a $k$-basis of $L$, then
      a preimage is given by the \emph{generalized trace form} $\Pmn\{L,x\}$, defined as
      \begin{equation}\label{eqspurform}
        \boxed{
          \Pmn\{L,x\}:=\text{Tr}_{L[X_i]/k[X_i]}
          \left[\frac{1}{x}\Bigl(e_1X_1+\ldots+e_nX_n\Bigr)^m\right]
          \in k[X_1,\ldots,X_n].
        }
      \end{equation}
  \end{enumerate}
\end{thm}

\vspace{\abstand}

\begin{proof}[Outline of proof:]
  For details refer to \cite{chris} or \cite{larsphd}. ---
  First of all, the canonical section $S_n\hookrightarrow\kranz$ of $p$ gives rise to a section of $p_*$
  which proves the surjectivity of $p_*$.\\[2 mm]
  Next we identify $S_n$ with $\AutN$ and as a consequence the short exact sequences
  \[
    \xymatrix@C=15mm@R=10mm{
      1 \ar[r] &
      {\mumn} \ar[r] \ar@{=}[d] &
      {\kranz} \ar[r]^{p} \ar@{=}[d] &
      {S_n} \ar[r] \ar@{=}[d] &
      1 \\
      1 \ar[r] &
      {\mumN} \ar[r] &
      {\kranzN} \ar[r]^{p} &
      {\AutN} \ar[r] &
      1, \\
    }
  \]
  twist the lower row with c, considered as a 1-cocycle from $G$ in $\kranzN$, and get
  the bijection (as explained in \ref{definonabel}\ref{defitwist})
  \[
    \Hc\left(G,\left(\mumN\right)_c\right)^{H^0\left(G,\left(\AutN\right)_c\right)}
    \stackrel{\sim}{\longrightarrow}
    p_*^{-1}[c].
  \]
  How to compute
  ${H^0\left(G,(\AutN)_c\right)}$
  and $\Hc\left(G,\left(\mumN\right)_c\right)$? --- First, it is not very difficult to see that
  \[
    \AutkL^{\text{opp}}\xrightarrow{\alpha}{H^0\left(G,(\AutN)_c\right)},\;\;\;
      a\mapsto[\varphi\mapsto\varphi\circ a]
  \]
  is a well defined isomorphism of groups.\\[3 mm]
  To compute $\Hc\left(G,\left(\mumN\right)_c\right)$, we take the short exact sequence
  \[
    1\rightarrow\mumN\rightarrow{\Kt}^N\xrightarrow{m}{\Kt}^N\rightarrow 1,
  \]
  twist it by $c$
  and consider the resulting long exact cohomology sequence:
  \begin{multline*}
    1\rightarrow
    H^0\left(G,\left(\mumN\right)_c\right)\rightarrow
    H^0\left(G,\KtcN\right)\xrightarrow{m}
    H^0\left(G,\KtcN\right)\xrightarrow{\delta} \\
    \Hc\left(G,\left(\mumN\right)_c\right)\rightarrow
    \Hc\left(G,\KtcN\right)\xrightarrow{m}\ldots
  \end{multline*}
  It can then be proved that
  \[
    \Lt\xrightarrow{\sim}H^0\left(G,\KtcN\right),\;\;\;
    x\mapsto(\varphi x)_{\varphi\in N}
  \]
  is a well defined group isomorphism
  and that the analogue of Hilbert 90 holds in this setting, i.e. we get the
  following commutative diagram with exact rows:
  \[
    \xymatrix@C=6mm{
      {H^0\left(G,\KtcN\right)} \ar[r]^{m} &
      {H^0\left(G,\KtcN\right)} \ar[r]^{\delta} &
      {\Hc\left(G,\left(\mumN\right)_c\right)} \ar[r] \ar@{=}[d] &
      {\Hc\left(G,\KtcN\right)} \ar@{=}[d] \\
      {\Lt} \ar[r]^{m} \ar[u]^{\wr} &
      {\Lt} \ar[r]^{\delta'} \ar[u]^{\wr} &
      {\Hc\left(G,\left(\mumN\right)_c\right)} \ar[r] &
      0. \\
    }
  \]
  This proves that we have an isomorphism
  ${\Lt}/{\Lt}^m\xrightarrow{\gamma}{\Hc\left(G,\left(\mumN\right)_c\right)}$,
  and a careful inspection of the involved maps shows that $\gamma$ is given by
  \[
    \bar{x}\mapsto
    \left(
      \left(
        \frac{s\sqrt[m]{s^{-1}\varphi x}}
          {\sqrt[m]{\varphi x}}
      \right)_{\varphi\in N}
    \right)_s
  \]
  and that we have a commutative diagram
  \[
    \xymatrix@R=3mm@C=8mm{
      \Hc(G,(\mumN)_c) \ar@{}[r]|{\times} &
      H^0(G,(\AutN)_c) & \ar[r] & &
      {\Hc(G,(\mumN)_c)} \\
      & & \\
      \LtLtm \ar@{}[r]|{\times} \ar[uu]^{\gamma}_{\wr} &
      {\text{Aut}_{k}(L)^{\text{opp}}} \ar[uu]^{\alpha}_{\wr} & \ar[r] & &
      \LtLtm \ar[uu]^{\gamma}_{\wr} \\
      {(\bar{x}} \ar@{}[r]|{,} & {a)} & \ar@{|->}[r] & & {\overline{a(x)}.}
    }
  \]
  Because a right-action of $\Aut_{k}(L)^{\text{opp}}$ is a left-action of $\Aut_{k}(L)$,
  this proves (i).\\[5 mm]
  To prove (ii), consider the $n\times n$-matrix
  $B:=\left(\frac{\varphi e_j}{\sqrt[m]{\varphi x}}\right)_{\varphi\in N,j\in\{1,\ldots,n\}}$.
  One can show that $B$ is regular and by a direct computation see that $\Pmn(BX)=\Pmn\{L,x\}$, i.e.
  that $B$ defines an isomorphism from $\Pmn\{L,x\}$ to $\Pmn$.
  So $\Pmn\{L,x\}$ is indeed a $\FKz/\Fkz$-form of $\Pmn$.
  Then we can use the isomorphism $B$ to compute $\vartheta(\Pmn\{L,x\})$ using the definition
  in \ref{satzdefitheta} and compare the result with what \eqref{eqcocycle} gives us to show
  $\vartheta(\Pmn\{L,x\})=(L,x)$,
  and this completes the proof of the theorem.
\end{proof}

\vspace{\abstand}

\begin{satz}\label{satzfrobenius2}
  Consider the particular case $k=\F_q$, and let $\bar{b}$ be an arbitrary class in $\Hc(G,\mumn)$,
  given by a pair $(L,x)$ in the sense of \ref{thmkranz}\ref{kranziso}.\\[1 mm]
  Then $L=\prod_{i=1}^rL_i$ for subfields $L_i$ of $K=\bar{\F}_q$ of finite degree
  $n_i$ over $\F_q$, and $x=(x_1,\ldots,x_r)$ with $x_i\in L_i^\times$.
  Choose $m$-th roots $y_i$ of $x_i$ in $K$.\\[1 mm]
  Then there is a 1-cocycle $b$ representing $\bar{b}$ that on the arithmetic Frobenius $f$
  is given by
  \[
    \fbox{$b_f=\prod_{i=1}^r\left[
      \left(y_i^{q^{n_i}-1},1,\ldots,1\right)
      \cdot
      z_{n_i}
    \right]$}
  \]
  (where $z_{n_i}\in S_{n_i}\leq S_n$ denotes the standard cycle of length $n_i$).
\end{satz}

\vspace{\abstand}

\begin{proof}[Proof:]
  Without loss of generality, we can assume $r=1$, i.e. $L=\F_{q^n}$. Then
  $N:=\homkLK=\{f^0,\ldots,f^{n-1}\}$ which we identify with $\{0,\ldots,n-1\}$,
  and it is clear that $f$ acts on $N$ as the standard cycle $z_n$.\\[2 mm]
  For $i\in\{0,\ldots,n-1\}$ we choose $y_1^{q^i}$ as an $m$-th root of $f^i(x)$
  and with these choices get
  \[
    \frac{f\sqrt[m]{f^{-1}\varphi x}}{\sqrt[m]{\varphi x}}=
    \left\{\begin{array}{ll}
      \displaystyle
      \frac{f\sqrt[m]{f^{-1}x}}{\sqrt[m]{x}}
      =\frac{\left(\sqrt[m]{f^{n-1}x}\right)^q}{y}
      =\frac{\left(y^{(q^{n-1})}\right)^q}{y}
      =y_1^{q^{n-1}-1} & \text{if $\varphi=f^0$,} \\
      \displaystyle
      \frac{f\sqrt[m]{f^{i-1}x}}{\sqrt[m]{f^ix}}
      =\frac{\left(y^{(q^{i-1})}\right)^q}{y^{(q^i)}}
      =\frac{y^{(q^i)}}{y^{(q^i)}}
      =1 & \text{otherwise.}
    \end{array}\right.
  \]
  This proves the proposition because of \eqref{eqcocycle}.
\end{proof}

\vspace{\abstand}

\begin{bsp}\label{bspintro}
  Take $m=3$, $n=6$, $k=\F_7$
  and choose $b:=(L,x)$ with $L:=\F_{2401}\times\F_{49}$ and $x:=(\beta^2,\alpha)$,
  where $\alpha$ \resp $\beta$ are generators of $\F_{49}^\times$ \resp $\F_{2401}^\times$
  satisfying $\alpha^2+5\alpha+5=0$ \resp $\beta^4+5\beta^3+4\beta^2+\beta+5=0$.\\[1 mm]
  We have $n_1=4$ and $n_2=2$, choose third roots $y_1$ of $\frac{1}{\beta}$
  and $y_2$ of $\frac{1}{\alpha^2}$ and calculate:
  \[
    \begin{array}{lcl}
      y_1^{7^4-1} & = &
        \left(\frac{1}{\beta}\right)^{\frac{7^4-1}{3}}
        =\left(\frac{1}{\beta}\right)^{800}
        =\beta^{2400-800}
        =\beta^{1600}
        =\left(\beta^{400}\right)^4
        =5^4
        =2,\\
      y_2^{7^2-1} & = &
        \left(\frac{1}{\alpha^2}\right)^{\frac{7^2-1}{3}}
        =\left(\frac{1}{\alpha^2}\right)^{16}
        =\alpha^{48-2\cdot 16}
        =\alpha^{16}
        =\left(\alpha^8\right)^2
        =5^2
        =4.
    \end{array}
  \]
  Then because of \ref{satzfrobenius2}, we have a 1-cocycle $b$ representing the class $(L,x)$
  that on $f$ is given by
  \[
    b_f=(2,1,1,1,4,1)\cdot[1234][56]\in\mu_3{\textstyle\int}S_6,
  \]
  and computing \eqref{eqspurform} explicitly shows that $P\{L,x\}$ is the polynomial from \eqref{eqbsp1}.
\end{bsp}

\vspace{\abstand}

\begin{cor}\label{corlangesequenz}
  Let $(L,x)$ represent an arbitrary class in $\Hc(G,\kranz)$ with
  $L=\prod_{i=1}^rL_i$ for subfields $L_i$ of $K$.
  Then we have the following short exact sequence of groups:
  \begin{equation}\label{eqautsequenz}
    \fbox{$
      1 \rightarrow
      {\displaystyle\prod_{i=1}^r(L_i\cap\mum)} \rightarrow
      \Autk\Bigl(\Pmn\{L,x\}\Bigr) \rightarrow
      \Bigl\{\Bigl.a\in\AutkL^{\text{opp}}\Bigr\vert\frac{ax}{x}\in{\Lt}^m\Bigr\} \rightarrow
      1
    $}
  \end{equation}
  In the special case $x=1\in\AutkL\backslash \LtLtm$ we get more precisely:
  \begin{equation}\label{eqautsemi}
    \fbox{$
      \displaystyle\Autk\Bigl(\Pmn\{L,x\}\Bigr)
      \cong
      \left(\prod_{i=1}^r(L_i\cap\mum)\right)\rtimes\AutkL^{\text{opp}}
    $}
  \end{equation}
  with the semi-direct product given by the following action:
  \begin{equation}\label{eqaction}
    \begin{array}{rclcc}
      \text{Aut}_k(L)^{\text{opp}} &
      \times &
      \left(\prod_{i=1}^r(L_i\cap\mum)\right) &
      \longrightarrow &
      \left(\prod_{i=1}^r(L_i\cap\mum)\right) \\[2 mm]
      {\scriptstyle (a} &
      {\scriptstyle ,} &
      {\scriptstyle x)} &
      {\scriptstyle\mapsto} &
      {\scriptstyle a^{-1}(x)}
    \end{array}
  \end{equation}
\end{cor}

\vspace{\abstand}

\begin{proof}[Proof:]
  We again set $N:=\homkLK$ and identify $N$ with $\{1,\ldots,n\}$, then
  take the short exact sequence
  \[
    1 \longrightarrow
    \mumN \longrightarrow
    \underbrace{\kranzN}_{=:H}\longrightarrow
    \AutN \longrightarrow
    1,
  \]
  twist by the 1-cocycle $b:=(L,x)$
  and look at the induced long exact cohomology sequence:
  \[
    1 \rightarrow
    H^0(G,(\mumN)_b) \rightarrow
    H^0(G,H_b) \rightarrow
    H^0(G,\AutN_b) \xrightarrow{\delta}
    \Hc(G,(\mumN)_b).
  \]
  Let $c$ denote the image of $b$ under $H\xrightarrow{p}\AutN\xrightarrow{s}H$
  (with $s$ the canonical section).
  Then it is easy to see that $\AutN_b=\AutN_c$ and $(\mumN)_b=(\mumN)_c$ and that
  we have a group isomorphism
  \[
    \eta:\prod_{i=1}^r(L_i\cap\mum)\xrightarrow{\sim}H^0(G,(\mumN)_c)),\;\;\;
    x\mapsto(\varphi x)_{\varphi\in N}.
  \]
  If we use this together with the isomorphism $\gamma$ defined in the proof of \ref{thmkranz}
  and apply \ref{corautform}, we get the following commutative diagram with exact rows:
  \[
    \xymatrix@R=13mm@C=4mm{
      1 \ar[r] &
      H^0(G,(\mumN)_c) \ar[r] \ar@{}[rd]|{=} &
      H^0(G,H_b) \ar[r] \ar@{}[rd]|{=} &
      H^0(G,\AutN_c) \ar[r]^{\delta} \ar@{}[rd]|{=} &
      \Hc(G,(\mumN)_c) \\
      1 \ar[r] &
      {\displaystyle\prod_{i=1}^r(L_i\cap\mum)} \ar[r] \ar[u]_{\eta}^{\wr} &
      \Autk\Bigl(\Pmn\{L,x\}\Bigr) \ar[r] \ar[u]^{\wr} &
      \Autk(L)^{\text{opp}} \ar@{-->}[r]_{\tilde{\delta}} \ar[u]_{\alpha}^{\wr} &
      \LtLtm \ar[u]_{\gamma}^{\wr} \\
    }
  \]
  The exactness of \eqref{eqautsequenz}
  now follows from the fact that $\tilde{\delta}$ is given by $a\mapsto\frac{ax}{x}$
  which can be proved by a direct computation using the explicit desciptions of
  $\alpha$, $\gamma$ and $\delta$.\\[3 mm]
  Consider now the case $x=1$. Then $b$ is induced from a 1-cocycle $c$ from $G$ in $S_n$,
  and it is easy to see that then $S_n\xrightarrow{s}\kranz$ defines a splitting
  of \eqref{eqautsequenz}.

  Finally, the explicit description of the action in \eqref{eqaction} is proved by comparing
  it with the natural action of $\AutN$ on $\mumN$ using the isomorphisms $\alpha$ and $\eta$.
\end{proof}

\vspace{\abstand}

\begin{bem}
  We deduced the exact sequences \eqref{eqautsequenz} and \eqref{eqautsemi}
  using our general formalism of coefficient extensions and the cohomological description of forms;
  as S{\l}adek and Weso{\l}owski show in \cite[Theorem 1.3.]{wes1} \resp \cite[Theorem 3.3.]{wes2},
  these results can also be obtained by a direct computation.
\end{bem}

\vspace{\abstand}

\begin{satz}\label{satzdefitildea}
  Define $\tA:=\tilde{A^m_n}:=\left(\kranz\right)/\mum$
  (where $\mum$ is embedded diagonally into $\kranz$).
  Then $\kranz\hookrightarrow\GLnK$ induces an embedding $\tA\hookrightarrow\PGLnK$ and an
  isomorphism onto $\Aut_{\FPZ}(\Pmn)$.
\end{satz}

\vspace{\abstand}

\begin{proof}[Proof:]
  The kernel of the canonical map $\kranz\hookrightarrow\GLnK\twoheadrightarrow\PGLnK$ is
  $\mum$, so that $\tA$ embeds into $\PGLnK$, and it is also clear that the image is contained
  in $\Aut_{\FPZ}(\Pmn)$. To prove equality, let $\bar{B}\in\PGLnK$ be an $\FPZ$-automorphism of $\Pmn$,
  represented by $B\in\GLnK$.
  By definition of morphisms in $\FPZ$, there exists $c\in \Kt$ with $\Pmn(BX)=c\cdot \Pmn$.
  Choose an $m$-th root $d$ of $c$ in $\Kt$ and put
  $B':=\frac{1}{d}B\in\GLnK$. Then
  $\Pmn(B'X)=\Pmn({\textstyle\frac{1}{d}}BX)=\left({\textstyle\frac{1}{d}}\right)^m\cdot \Pmn(BX)=\Pmn$,
  i.e. $B'$ is an $\FKz$-automorphism of $\Pmn$. Then $B'\in\kranz$ because of \ref{satzshioda}, and
  this proves the proposition because of $\bar{B'}=\bar{B}\in\PGLnK$.
\end{proof}

\vspace{\abstand}

\vspace{\abstand}

\begin{defi}
  Let $L$ be a commutative, separable $k$-algebra of degree $n$. We introduce the following
  equivalence relation $\sim_k$ on the set $\AutkL\backslash(\LtLtm)$:
  \[
    [x]\sim_k[x']:\Leftrightarrow
    \exists\lambda\in k^\times:
    \exists a\in\AutkL:
    \exists y\in \Lt:
    x'=\lambda\cdot a(x)\cdot y^m.
  \]
\end{defi}

\begin{thm}\label{thmclassification}\mbox{}\\[-6mm]
  \begin{enumerate}
    \item\label{isoim}
      Let $p$ denote the canonical projection $\kranz\twoheadrightarrow\tA$. Then $\vartheta$ induces
      an isomorphism
      \[
        \fbox{$E(\FPZ/\FPz,\Pmn)
          \stackrel{\sim}{\longrightarrow}
          \text{Im}\,\left(\Hc(G,\kranz)\xrightarrow{p_*}\Hc(G,\tA)\right)$}
      \]
    \item\label{simk}
      Two classes $(L,x)$ and $(L',x')$ in $\Hc(G,\kranz)$ have the same image under $p_*$ iff
      there exists $\lambda\in k^\times$ such that $(L,\lambda x)=(L',x')$.
      Because of \ref{thmkranz}\ref{kranziso}, this is equivalent to the existence of a $k$-isomorphism
      $\psi:L\xrightarrow{\sim}L'$ satisfying $x'\sim_k\psi(x)$.
    \item\label{efpz}
      $\vartheta$ induces an isomorphism of pointed sets
      \[
        \fbox{$
          E(\FPZ/\FPz,\Pmn)\stackrel{\sim}{\longrightarrow}
          \displaystyle\coprod_{[L]}\Bigl[\AutkL\backslash(\LtLtm)\Bigr]/\sim_k
        $}
      \]
      where $[L]$ runs over the isomorphism classes of commutative, separable $k$-algebras of degree $n$.
      The inverse map sends a pair $(L,x)$ to the isomorphism class of $\Pmn\{(L,x)\}$.
  \end{enumerate}
\end{thm}

\vspace{\abstand}

\begin{proof}
  The obvious morphism of coefficient extensions
  $(\Fkz\rightarrow\FKz)\rightarrow(\FPz\rightarrow\FPZ)$
  induces because of \ref{satzMKE} the commutative diagram of pointed sets
  \[
    \xymatrix@R=15mm@C=25mm{
      {E(\FKz/\Fkz,\Pmn)}
        \ar[rr]^{\textstyle\varphi}_{\textstyle[Q]\mapsto[Q]}
        \ar[d]_{\textstyle\vartheta'}^{\textstyle\wr} \ar@{}[rrd]|{=}
      & & {E(\FPZ/\FPz,\Pmn)}
        \ar@{^{(}->}[d]^{\textstyle\vartheta} \\
      {\Hc(G,\kranz)}
        \ar[rr]_{\textstyle p_*}^{\textstyle\left(a_s\right)
        \mapsto\left(\bar{a}_s\right)}
      & & {\Hc(G,\tA)}.
    }
  \]
  This immediately implies
  $\text{Im}\,\vartheta\supseteq\text{Im}\,p_*$,
  and it also shows that in order to prove
  $\text{Im}\,\vartheta\subseteq\text{Im}\,p_*$, it suffices to see that $\varphi$ is surjective.
  So let $Q$ be an arbitrary $\FPZ/\FPz$-form of $\Pmn$,
  so that there are $A\in\GLnK$ and $\lambda\in\Kt$ with $\Pmn(AX)=\lambda\cdot Q$.
  Choose an $m$-th root $\mu$ of $\lambda$ and put $A':=(1/\mu)\cdot A$, then
  $\Pmn(A'X)=(1/\mu)^m\cdot\Pmn(AX)=(1/\lambda)\cdot\lambda\cdot Q=Q$,
  which means that $A'$ defines an isomorphism from $Q$ to $\Pmn$ in $\FKz$, i.e.
  $Q$ is also a $\FKz/\Fkz$-form of $\Pmn$, and $[Q]$ is a preimage of $[Q]$ under $\varphi$.
  This completes the proof of \ref{isoim}.

  Now let $b=(L,x)$ be an arbitrary element of $\Zc(G,\kranz)$.
  According to \ref{definonabel}\ref{defitwist}, those classes in $\Hc(G,\kranz)$ that are mapped
  to the same class as $[b]$ under $p_*$ are exactly those of the form
  $(a_s\cdot b_s)$ for $(a_s)\in\Hc(G,\mu_{m,b})$.

  It is easy to see that $\mu_{m,b}=\mu_m$, because $\mu_m$ lies in the center of $\kranz$.
  Therefore we have the well-known Kummer isomorphism
  $k^\times/(k^\times)^m\stackrel{\sim}{\longrightarrow}\Hc(G,\mu_{m,b})$,
  $\lambda\mapsto\left({s\sqrt[m]{\lambda}}/{\sqrt[m]{\lambda}}\right)_s$,
  and a straightforward calculation using \eqref{eqcocycle} shows
  $\left(\frac{s\sqrt[m]{\lambda}}{\sqrt[m]{\lambda}}\cdot b_s\right)=(L,\lambda x)$ in $\Hc(G,\kranz)$
  which proves \ref{simk}.

  Finally, \ref{efpz} is a direct consequence of \ref{isoim}, \ref{simk} and \ref{thmkranz}.
\end{proof}

\vspace{\abstand}

%%%%%%%%%%%%%%%%%%%%%%%%%%%%%%%%%%%%%%%%%%%%%%%%%%%%%%%%%%%%%%%%%%%%%%%%%%%%%%%%%%%%%%%%%%%%%%%%%%%%%%%%%%%%%%
%%%%%%%%%%%%%%%%%%%%%%%%%%%%%%%%%%%%%%%%%%%%%%%%%%%%%%%%%%%%%%%%%%%%%%%%%%%%%%%%%%%%%%%%%%%%%%%%%%%%%%%%%%%%%%
%%%%%%%%%%%%%%%%%%%%%%%%%%%%%%%%%%%%%%%%%%%%%%%%%%%%%%%%%%%%%%%%%%%%%%%%%%%%%%%%%%%%%%%%%%%%%%%%%%%%%%%%%%%%%%
%%%%%%%%%%%%%%%%%%%%%%%%%%%%%%%%%%%%%%%%%%%%%%%%%%%%%%%%%%%%%%%%%%%%%%%%%%%%%%%%%%%%%%%%%%%%%%%%%%%%%%%%%%%%%%
%%%%%%%%%%%%%%%%%%%%%%%%%%%%%%%%%%%%%%%%%%%%%%%%%%%%%%%%%%%%%%%%%%%%%%%%%%%%%%%%%%%%%%%%%%%%%%%%%%%%%%%%%%%%%%
%%%%%%%%%%%%%%%%%%%%%%%%%%%%%%%%%%%%%%%%%%%%%%%%%%%%%%%%%%%%%%%%%%%%%%%%%%%%%%%%%%%%%%%%%%%%%%%%%%%%%%%%%%%%%%
%%%%%%%%%%%%%%%%%%%%%%%%%%%%%%%%%%%%%%%%%%%%%%%%%%%%%%%%%%%%%%%%%%%%%%%%%%%%%%%%%%%%%%%%%%%%%%%%%%%%%%%%%%%%%%
%%%%%%%%%%%%%%%%%%%%%%%%%%%%%%%%%%%%%%%%%%%%%%%%%%%%%%%%%%%%%%%%%%%%%%%%%%%%%%%%%%%%%%%%%%%%%%%%%%%%%%%%%%%%%%
%%%%%%%%%%%%%%%%%%%%%%%%%%%%%%%%%%%%%%%%%%%%%%%%%%%%%%%%%%%%%%%%%%%%%%%%%%%%%%%%%%%%%%%%%%%%%%%%%%%%%%%%%%%%%%

\section{$l$-adic cohomology of the Fermat hypersurface}

\vspace{\abstand}

\begin{satz}\label{cormot}
  Let $S$ be a group,
  $A$ a finite abelian group of order $N$ and exponent $m$,
  $\Z[\frac{1}{N},\mum]\subseteq R\subseteq\C$ a ring,
  $\Mot$ a pseudoabelian, $R$-linear category,
  and $M\in\text{Ob}(\Mot)$ an object with a \emph{right}-$(A\rtimes S)$-action.
  Then an easy calculation, using basic representation theory of finite groups 
  (see for example \cite[pp675]{lang}) shows:
  \begin{enumerate}
    \item
      The set $\{p_\chi\,\vert\,\chi\in\check{A}\}$ with
      $p_\chi:=\frac{1}{n}\sum_{a\in A}\chi(a)^{-1}a\in R[A]$
      forms a complete system of idempotents in $R[A]$.
    \item
      If we consider the $p_\chi$ via $A=A^{\text{opp}}\rightarrow\Aut(M)$ as projectors on $M$
      and denote the parts cut out by $M_\chi$, then
      $M=\bigoplus_{\chi\in\check{A}}M_\chi$.
    \item\label{cormotdiagrams}
      For $s\in S$, $b\in A$ and $\chi\in\check{A}$, we have the following commutative diagrams:
      \[
      \xymatrix@C=2.2cm@R=1.5cm{
      M_{\op{s}{\chi}} \ar[r]^{p_\chi\cdot s\cdot p_{\op{s}{\chi}}} \ar[d] \ar@{}[rd]|{=} & M_\chi \ar[d] &
        M_\chi \ar[r]^{\chi(b)} \ar[d] \ar@{}[rd]|{=} & M_\chi \ar[d]\\
      M \ar[r]_{s} & M & M \ar[r]_{b} & M
      }
      \]
    \item
      For $\chi\in\check{A}$, let $[\chi]\subseteq\check{A}$ be the orbit of $\chi$
      under the action of $S$, and set $M_{[\chi]}:=\bigoplus_{\chi\in[\chi]}M_{\chi}$.
      Then $A\rtimes S$ acts on $M_{[\chi]}$, and if $\check{A}=[\chi_1]\sqcup\ldots\sqcup[\chi_r]$,
      then $M=M_{[\chi_1]}\oplus\ldots\oplus M_{[\chi_r]}$
      is a decomposition of $M$ into $(A\rtimes S)$-invariant subobjects.
  \end{enumerate}
\end{satz}

\vspace{\abstand}

\begin{bspe}\label{bspzerl}
  Let $X/k$ be a smooth, projective variety with a left-$(A\rtimes S)$-action.
  \begin{enumerate}
    \item
      Take for $\Mot$ the category $\MotEk$ of Grothendieck-motives over $k$ with coefficients in $E$, and
      let $h(X)$ be the motive of $X$. Then $A\rtimes S$ acts on $h(X)$ from the right, and we get
      \[
        h(X)=\bigoplus_{\chi\in\check{A}}h(X)_\chi.
      \]
    \item
      Let $l\neq\text{char}(k)$ be a prime with $l\equiv 1\pmod{m}$ and $i\in\N_0$.
      Then $\Q_l$ contains $\mum$ and, we can choose $\Q(\mum)\hookrightarrow\Q_l$ by which
      $\Mot:=\Qlk$ becomes a \mbox{(pseudo-)}abelian, $\Q(\mum)$-linear category.
      If $V\in\text{Ob}(\Mot)$ denotes the $\Q_l$-$G_k$-representation $\Het(X\times_k\bar{k},\Q_l)$,
      we get a right-$(A\rtimes S)$-action on $V$ and therefore a decomposition of $V$ in
      $\Q_l$-$G_k$-representations $V_\chi$, explicitly given by
      \begin{equation}\label{eqvchi}
        V_\chi=\{v\in V\,|\,\forall a\in A:\,v\cdot a=\chi(a)\cdot v\}.
      \end{equation}
  \end{enumerate}
\end{bspe}

\vspace{\abstand}

\noindent
Now consider the case $n\geq 2$, $k=\F_q$, 
$q\equiv 1\pmod{m}$ ($\Longleftrightarrow\mum\subseteq\F_q^\times$),
$p:=\text{char}\,(k)>m$,
and choose a prime $l\neq p$, $l\equiv 1\pmod{m}$.

\vspace{\abstand}

\begin{lemmadefi}\label{lemmadefiamn}
  Let ${\mathcal X}:={\mathcal X}^m_n$ be the \emph{Fermat hypersurface} associated to $\Pmn$,
  $A$ the finite abelian group $\mumn/\mum$ (so that $\tilde{A}=A\rtimes S_n$),
  $\zeta:=\mbox{e}^{\frac{2\pi i}{m}}\in\mum\subset\C$,
  and choose an embedding $\mum\hookrightarrow\Q(\zeta)$.
  Then
  \[
    \check{A}\cong
    \left\{
      \fata=(a_1,\ldots,a_n)\in(\Z/m\Z)^n\ \left|\ \sum_{i=1}^na_i=0\in\Z/m\Z\right.
    \right\}
  \]
  via
  \begin{equation}\label{eqdual}
    (\fata\def\fata{\mbox{\boldmath $a$}},\underbrace{[(\zeta_1,\ldots,\zeta_n)]}_{\in A})
    \mapsto
    \prod_{i=1}^n\zeta_i^{a_i}.
  \end{equation}
  Define
  \[
    A^m_n:=\left\{
      \fata\in\check{A}\ \left|\ \forall i\in\{1,\ldots,n\}\ :\ a_i\neq 0\in\Z/m\Z\right.
    \right\}.
  \]
  If we set $\bar{\mathcal X}:={\mathcal X}\times_k\bar{k}$ and
  $V:=V^m_n:={\Hetn(\bar{\mathcal X},\Q_l)}$, according to \ref{bspzerl}(ii) we get a canonical decomposition
  of the $\Q_l$-$G_k$-representation $V$ as
  \[
    V=\bigoplus_{\chi\in\check{A}}V_\chi.
  \]
\end{lemmadefi}

\vspace{\abstand}

\begin{defi}\label{defijacobisumme}
  Choose a character $\chi$ of $k^\times$ of exact order $m$,
  and for $\fata\in A^m_n$ define the {\em Jacobi sum of dimension $(n-2)$ and degree $m$ of $\fata$}
  as
  \begin{multline*}
    {\mathcal J}_k^m(\fata):={\mathcal J}(\fata)\\
    :=
    (-1)^n\sum_{\begin{array}{c}
      \scriptstyle (v_2,\ldots,v_n)\in k^\times\times\ldots\times k^\times\\
      \scriptstyle 1+v_2+\ldots+v_n=0
    \end{array}}
    \chi(v_2)^{a_2}\cdot\chi(v_3)^{a_3}\cdot\ldots\cdot\chi(v_n)^{a_n}\in\Q(\zeta)\subseteq\Q_l.
  \end{multline*}
\end{defi}

\vspace{\abstand}

\begin{satz}\label{satzjacobi}(Weil \cite{weil})\\
  If $\sigma\in S_n$ is a permutation and $\fata\in A^m_n$, then
  ${\mathcal J}(\sigma\fata)={\mathcal J}(\fata)$,
  i.e. ${\mathcal J}(\fata)$ only depends on $[\fata]\in S_n\backslash A^m_n$.
\end{satz}

\vspace{\abstand}

\begin{thm}\label{satzkohom}(Deligne \cite[I., \S 7]{deligne})\\
  \[
    V=\underbrace{\left(\bigoplus_{a\in A^m_n}V_{\fata}\right)}_{=:V^m_{n,\text{prim}}=:V_{\text{prim}}}\oplus
    \left\{\begin{array}{cl}
      0 & \text{, for $n$ odd,}\\
      \Q_l(-\frac{n-2}{2}) & \text{, for $n$ even,}
    \end{array}\right.
  \]
  $\dim_{\Q_l}V_{\fata}=1$,
  and the geometric Frobenius $\FX$ acts on $V_{\fata}$ as multiplication by ${\mathcal J}(\fata)$.
\end{thm}

\vspace{\abstand}

\begin{bsp}\label{bspjacobi}
  For $m=3$, $n=6$, $k=\F_7$ as in \ref{bspintro}, we have
  $V=V_{\{0,0,0,0,0,0\}}\oplus V_{\text{prim}}$,
  $V_{\text{prim}}=V_{\{1,1,1,1,1,1\}}\oplus V_{[1,1,1,2,2,2]}\oplus V_{\{2,2,2,2,2,2\}}$,
  and if we define $\chi:\F_7^\times\rightarrow\C$ by $\chi(3):=\zeta$, then
  \[
    \begin{array}{lll}
      {\mathcal J}(1,1,1,1,1,1) & = & -56-21\zeta, \\
      {\mathcal J}(1,1,1,2,2,2) & = & 49, \\
      {\mathcal J}(2,2,2,2,2,2) & = & -35+21\zeta. \\
    \end{array}
  \]
\end{bsp}

\vspace{\abstand}

\begin{lemma}\label{lemmahyperebene}
$\tilde{A}$ acts trivially on $V_{(0,\ldots,0)}$.
\end{lemma}

\vspace{\abstand}

\begin{proof}
  If $n$ is odd, $V_{(0,\ldots,0)}=0$ according to \ref{satzkohom}, and there is nothing to prove.
  So let $n$ be even.\\[2 mm]
  In the case $n=2$, we have $\dim{\mathcal X}=0$, i.e. $V$ is a $\Q_l$-algebra
  containing a subring $R$ that is isomorphic to $\Q_l$ and on which all automorphisms of $\mathcal X$
  act trivially. So in particular $A$ acts trivially on $R$ which means $R\subseteq V_{(0,0)}$. 
  But $\dim_{\Q_l}V_{(0,0)}=1$ according to \ref{satzkohom}, so $R=V_{(0,0)}$ follows.\\[2 mm]
  In the case $n\geq 4$, 
  let $[H]$ denote the class of the smooth hyperplane section $\{x_n=0\}$ in $\text{CH}^1_{\Q}({\mathcal X})$,
  and let $\gamma\in{{\text{H}_{\text{\'{e}t}}^2}(\bar{\mathcal X},\Q_l(1))}$ 
  denote the corresponding class in cohomology.
  Hard Lefschetz tells us that multiplication by $\gamma^{n-2}$ induces an isomorphism
  $\Hetnull\xrightarrow{\sim}{{\text{H}_{\text{\'{e}t}}^{2n-4}}(\bar{\mathcal X},\Q_l(n-2))}$
  which in particular implies $\left[\gamma(-1)\right]^{\frac{n-2}{2}}\neq 0\in V$.
  Obviously, for all automorphisms $s\in\tilde{A}$ we have $s^*[H]=H$, 
  so that $\tilde{A}$ (and thus its subgroup $A$) acts trivially one the one-dimensional subspace
  $\left\langle\left[\gamma(-1)\right]^{\frac{n-2}{2}}\right\rangle$ of $V$.
  Then \ref{satzkohom} implies that this subspace is equal to $V_{(0,\ldots,0)}$, which completes the
  proof of the lemma.
\end{proof}

\vspace{\abstand}

\noindent
The following proposition is surprisingly hard to prove and constitutes the technical heart of this paper:

\begin{satz}\label{lemmaminuseins}
  Let $\tau\in S_n$ be a \emph{transposition} and $\fata\in A^m_n$ with $^\tau\fata=\fata$,
  so that $\tau$ induces a $\Q_l$-linear involution on $V_{\fata}$ according to 
  \ref{cormot}\ref{cormotdiagrams}. 
  This involution is multiplication by $(-1)$.
\end{satz}

\vspace{\abstand}

\begin{proof}
  Without loss of generality, we may assume $\tau=[12]$ (and therefore $a_1=a_2$). We prove the proposition
  seperately for the three cases $n=2$, $n=3$ and $n\geq 4$:
  \begin{itemize}
    \item {\em\underline{$n=2$}:}
      If $m$ is \emph{odd}, we have $A^m_2=\emptyset$, and nothing is to prove. So let $m$ be \emph{even}
      which in particular implies $\fata=(\frac{m}{2},\frac{m}{2})$. 
      Because of \ref{cormot}\ref{cormotdiagrams}, we have the decomposition
      \[
        V=V_{(0,0)}\oplus V_{(\frac{m}{2},\frac{m}{2})}\oplus
        \underbrace{\bigoplus_{j\in\{1,\ldots,m-1\}\setminus\{\frac{m}{2}\}}
        V_{(j,m-j)}}_{=:V'}
      \]
      of $V$ into $\Q_l$-$G_k$-subspaces.
      Lemma \ref{lemmahyperebene} tells us that $\tau^*$ has trace 1 on $V_{(0,0)}$,
      and \ref{cormot}\ref{cormotdiagrams} shows that $\tau^*$ has trace 0 on $V'$.
      Finally, a simple calculation shows that the automorphism $\tau$ of $\bar{\mathcal X}$ 
      has no fixed points, and therefore Lefschetz trace formula gives us
      \[
        0=\text{Tr}\,(\tau^*)=1+\text{Tr}\,(\tau^*|V_{(\frac{m}{2},\frac{m}{2})})+0
        \;\;\Longrightarrow\;\;
        \text{Tr}\,(\tau^*|V_{(\frac{m}{2},\frac{m}{2})})=-1.
      \]
      Seeing as $V_{(\frac{m}{2},\frac{m}{2})}$ is one-dimensional because of \ref{satzkohom}, 
      this completes the proof for $n=2$.
    \item {\em\underline{$n=3$}:}
      Again, we first compute the number $(\Delta_{\bar{\mathcal X}},\Gamma_\tau)$ of fixed points of $\tau$
      by a simple calculation:
      \begin{equation}\label{eqdelta}
        (\Delta_{\bar{\mathcal X}},\Gamma_\tau)=\left\{
          \begin{array}{cl}
            m & \text{, for $m$ even,}\\
            m+1 & \text{, for $m$ odd.}
          \end{array}
        \right.
      \end{equation}
      Lefschetz trace formula implies
      \begin{equation}\label{eq9}
        (\Delta_{\bar{\mathcal X}},\Gamma_\tau)=
        \underbrace{\text{Tr}\,(\tau^*|\Hetnull)}_{=:t_0}-
        \underbrace{\text{Tr}\,(\tau^*|\Heteins)}_{=:t_1}+
        \underbrace{\text{Tr}\,(\tau^*|\Hetzwei)}_{=:t_2}.
      \end{equation}
      Because of $n>2$, 
      $\bar{\mathcal X}$ is irreducible, 
      i.e. $\Hetnull$ is isomorphic to the trivial $\Q_l$-$G_k$-representation
      $\Q_l$, which implies $t_0=1$ and then $t_2=1$ by Poincar\'{e} duality.
      Therefore \eqref{eq9} gives us:
      \begin{equation}\label{eq10}
        t_1=2-(\Delta_{\bar{\mathcal X}},\Gamma_\tau).
      \end{equation}
      Now put $N:=\bigl\{\bigl.(a_1,a_2,a_3)\in A^m_3\,\bigr\vert\,a_1=a_2\bigr\}$ and note
      \begin{equation}\label{eqnumn}
        \#N=\left\{
          \begin{array}{cl}
            m-2 & \text{, for $m$ even,}\\
            m-1 & \text{, for $m$ odd.}
          \end{array}
        \right.
      \end{equation}
      It is
      \[
        t_1=
        \text{Tr}\,\left(\tau^*\left\vert\bigoplus_{a\in A^m_3}V_{\fata}\right.\right)=
        \sum_{\fata\in A^m_n}\underbrace{\text{Tr}\,(\tau^*|V_{\fata})}_{\in\{-1,1\}}
          +\underbrace{\text{Tr}\,\left(\tau^*\,\left\vert\,
            \bigoplus_{\fata\in A^m_3\setminus N}V_{\fata}\right.\right)}
          _{\stackrel{\ref{cormot}\ref{cormotdiagrams}}{=}0}
      \]
      in $[-\#N,\#N]$, and $t_1=-\#N$ iff the lemma is true for all $\fata\in N$.
      We check:
      \[
        t_1
        \stackrel{\eqref{eq10}}{=}
        2-(\Delta_{\bar{\mathcal X}},\Gamma_\tau)
        \stackrel{\eqref{eqdelta}}{=}
        2-\left\{\begin{array}{c}m\\m+1\end{array}\right\}
        =
        \left\{\begin{array}{c}2-m\\1-m\end{array}\right\}
        \stackrel{\eqref{eqnumn}}{=}
        -\#N.
      \]
    \item {\em\underline{$n\geq 4$}:}
      Let $r,s\in\N^{\geq 3}$ with $r+s=n+2$. We use Shioda's ``inductive structure of
      Fermat varieties'' (\cite[p.179]{shioda79}, \cite{shioda82}):
      According to Shioda, we have the following commutative diagram:\\
      \[\xymatrix@C=5mm{
        {\beta^{-1}(Y)} \ar@{^{(}->}[r]^{j'} \ar[d]_{\beta'} \ar@{}[rd]|{\Box} &
          {Z_{r,s}^m} \ar[r]^-{\pi} \ar[d]_{\beta} \ar[rd]^{\psi} &
          {Z_{r,s}^m/\mum} \ar[d]^{\bar{\psi}} \\
        {\underbrace{{\mathcal X}_{r-1,K}^m\times{\mathcal X}_{s-1,K}^m}_{=:Y}} \ar@{^{(}->}[r]_-{j} &
          {{\mathcal X}_{r,K}^m\times{\mathcal X}_{s,K}^m} \ar@{..>}[r]_-{\varphi} &
          {\bar{\mathcal X}} &
          {{\mathcal X}_{r-1,K}^m\coprod{\mathcal X}_{s-1,K}^m,} \ar@{_{(}->}[l]^-{i}
      }\]
      \mbox{}\\
      where the maps are defined as follows:\\[2 mm]
      \begin{tabular}{ccl}
        $\varphi$ & : & rational map, defined by \\
        & & $([x_1:\ldots:x_r],[y_1:\ldots:y_s])\mapsto$ \\
        & & \hspace{1 cm}$[x_1y_s:\ldots:x_{r-1}y_s:\varepsilon x_ry_1:\ldots:\varepsilon x_ry_{s-1}]$
          (with $\varepsilon^m=-1$),\\
        $j$ & : & $([x_1:\ldots:x_{r-1}],
          [y_1:\ldots:y_{s-1}])\mapsto$\\
        & & \hspace{1 cm}$([x_1:\ldots:x_{r-1},0],[y_1:\ldots:y_{s-1},0])$,\\
        $i$ & = & $i_1\coprod i_2:
          \left\{\begin{array}{l}
            i_1([x_1:\ldots:x_{r-1}])=[x_1:\ldots:x_{r-1}:0:\ldots:0]\\
            i_2([y_1:\ldots:y_{s-1}])=[0:\ldots:0:y_1:\ldots:y_{s-1}]
          \end{array}\right.$,\\
        $\beta$ & : & blowing up of ${\mathcal X}_{r,K}^m\times{\mathcal X}_{s,K}^m$ along $Y$,\\
        $\beta'$ & : & restriction of $\beta$,\\
        $j'$ & : & embedding,\\
        $\pi$ & : & quotient morphism,\\
        $\bar{\psi}$ & : & blowing up of $\bar{\mathcal X}$ 
        along ${\mathcal X}_{r-1,K}^m\coprod{\mathcal X}_{s-1,K}^m$,\\
        $\psi$ & = & $\varphi\circ\beta=\bar{\psi}\circ\pi$.
      \end{tabular}
      \mbox{}\\[6 mm]
      Shioda defines 
      $A^m_{r,s}:=\{(\fatb,\fatc)\in A^m_r\times A^m_s\,|\,b_r+c_s=0\}$,
      considers the maps
      \[
        \begin{array}{lcllcl}
          A^m_{r,s} & \stackrel{\sharp}{\longrightarrow} & A^m_n, &
          \fatb\,\sharp\,\fatc & := & (b_1,\ldots,b_{r-1},c_1,\ldots,c_{s-1})\ \text{and} \\
          A^m_{r-1}\times A^m_{s-1} & \stackrel{*}{\longrightarrow} & A^m_n, &
          \fatbp\star\fatcp & := & (b'_1,\ldots,b'_{r-1},c'_1,\ldots,c'_{s-1}),
        \end{array}
      \]
      shows that 
      $A^m_{r,s}\coprod\left(A^m_{r-1}\times A^m_{s-1}\right)\xrightarrow{\sharp\,\sqcup\,*}A^m_n$
      is a bijection and proves that
      for $(\fatb,\fatc)\in A^m_{r,s}$ and $(\fatbp,\fatcp)\in A^m_{r-1}\times A^m_{s-1}$,
      we get isomorphisms
      \[\xymatrix@C=1.0cm{
        V_{\fatb}\otimes V_{\fatc}
        \ar[r]^-{\psi_*\beta^*}_-{\sim} &
        V_{\fatb\,\textstyle\sharp\,\fatc} &
        *\txt{and} &
        V_{\fatbp}\otimes V_{\fatcp}(-1)
        \ar[r]^-{\psi_*{j'}_*{\beta'}^*}_-{\sim} &
        V_{\fatbp\textstyle\star\fatcp}.
      }\]
      In particular, we can choose $r:=3$ and $s:=n-1$, and
      look at the following diagram:\\
      \xymatrix@C=4.5mm{
        & {\beta^{-1}(Y)} \ar@{^{(}->}[rr]^{j'} \ar[d]_{\beta'} & &
          {Z_{3,n-1}^m} \ar[d]_{\beta} \ar[rrrd]^{\psi} \\
        & {Y} \ar@{^{(}->}[rr]_j|!{[rd];[rru]}\hole & &
          {{\mathcal X}_{3,K}^m\times{\mathcal X}_{n-1,K}^m} \ar@{..>}[rrr]_{\varphi} & & & 
          {\bar{\mathcal X}} \\
        {\beta^{-1}(Y)} \ar@{^{(}->}[rr]^{j'} \ar[d]_{\beta'} \ar[uur]^(.65){\tau'} & &
          {Z_{3,n-1}^m} \ar[d]_{\beta} \ar[rrd]^{\psi}  \ar[uur]^(.65){\tau'} \\
        {Y} \ar@{^{(}->}[rr]_j \ar[uur]^(.65){\tau\times 1}|!{[u];[urr]}\hole & &
          {{\mathcal X}_{3,K}^m\times{\mathcal X}_{n-1,K}^m} \ar@{..>}[rr]_{\varphi}
          \ar[uur]^(.65){\tau\times 1}|!{[u];[rr]}\hole & &
          {\bar{\mathcal X}} \ar[uurr]_(.65){\tau}
      }
      \mbox{}\\
      The morphism $\tau\times 1$ maps $Y$ isomorphically to $Y$, 
      so that we get an induced $\tau'$ 
      on the blowing up that commutes with $\tau\times 1$.
      Furthermore, one immediately sees
      $(\tau\times 1)\circ j=j\circ(\tau\times 1)$ and
      $\tau\circ\varphi=\varphi\circ(\tau\times 1)$ which shows that the diagram commutes.\\[2 mm]
      Because $(\sharp\,\sqcup*)$ is a bijection, $\fata$ is either of the form
      $\fatb\sharp\fatc$ or of the form $\fatbp\star\fatcp$,
      and our assumption $\op{\tau}\fata=\fata$ implies
      $b_1=b_2$ \resp $b'_1=b'_2$. 
      Because $\tau$ is an involution, we have $\tilde{\tau}_*=\tilde{\tau}^*$,
      and we get\\
      \hspace*{5 mm}$\tau^*(\psi_*\beta^*)=(\tau_*\psi_*)\beta^*=(\psi_*{\tau'}_*)\beta^*
        =\psi_*({\tau'}^*\beta^*)=(\psi_*\beta^*)(\tau\times 1)^*$
      \resp\\
      \hspace*{5 mm}$\tau^*(\psi_*j'_*{\beta'}^*)=(\tau_*\psi_*j'_*){\beta'}^*
        =(\psi_*j'_*{\tau'}_*){\beta'}^*
        =\psi_*j'_*({\tau'}^*{\beta'}^*)$
      \begin{flushright}
        $=(\psi_*j'_*{\beta'}^*)(\tau\times 1)^*$,
      \end{flushright}
      which shows that the following diagrams commute:
      \[\xymatrix@C=0.8cm@R=0.55cm{
        {V_{\fatb}\otimes V_{\fatc}}
          \ar[r]^-{\psi_*\beta^*}_-{\sim}
          \ar[dd]_{(\tau\times 1)^*}^{\wr} &
          {V_{\fatb\,\textstyle\sharp\,\fatc}}
          \ar[dd]^{\tau^*}_{\wr} & &
        {V_{\fatbp}\otimes V_{\fatcp}(-1)}
          \ar[r]^-{\psi_*{j'}_*{\beta'}^*}_-{\sim}
          \ar[dd]_{(\tau\times 1)^*}^{\wr} &
          {V_{\fatbp\textstyle\star\fatcp}}
          \ar[dd]^{\tau^*}_{\wr} \\
        & & {\text{\resp}} \\
        {V_{\fatb}\otimes V_{\fatc}}
          \ar[r]^-{\psi_*\beta^*}_-{\sim} &
          {V_{\fatb\,\textstyle\sharp\,\fatc}} & &
        {V_{\fatbp}\otimes V_{\fatcp}(-1)}
          \ar[r]^-{\psi_*{j'}_*{\beta'}^*}_-{\sim} &
          {V_{\fatbp\textstyle\star\fatcp}}
      }\]
      As we already treated the cases ``$n=2$'' and ``$n=3$'',
      we know that $(\tau\times 1)^*$ is multiplication by $(-1)$, 
      and the claim follows from the commutativity of the respective diagram.
  \end{itemize}
  So the claim is proved for all $n\geq 2$.
\end{proof}

\vspace{\abstand}

\begin{cor}\label{corbasis}
  Let $\fata\in A^m_n$, and fix $v\in V_{\fata}\setminus\{0\}$.
  For each $\fatb\in[\fata]:=S_n\cdot\fata\in S_n\backslash A^m_n$, 
  choose $\sigma_b\in S_n$ with $\sigma_b^{-1}(\fata)=\fatb$
  and set $v_b:=\sign\,(\sigma_b)\cdot\sigma_b^*v$.\\[2 mm]
  Then $\{v_b\}$ is a basis of $V_{[\fata]}:=\bigoplus_{\fatb\in[\fata]}V_{\fatb}$, 
  independant of the chosen $\sigma_b$, such that for all
  $s:=(\zeta_1,\ldots,\zeta_n)\cdot\sigma\in {\tilde{A}}$ and all $\fatb\in[\fata]$:
  \begin{equation}\label{eqvb}
    \fbox{$
    s^*v_b=(\zeta_1^{b_{\sigma(1)}}\cdot\ldots\cdot\zeta_n^{b_{\sigma(n)}})\cdot\sign\,(\sigma)
      \cdot v_{\sigma^{-1}b}
    $}.
  \end{equation}
\end{cor}

\vspace{\abstand}

\begin{proof}[Proof:]
  First of all, $\forall\fatb\in[\fata]:v_b\in V_{\fatb}\setminus\{0\}$, 
  because $v\neq 0$ and $\sigma^*$ is an
  automorphism of $V_{[\fata]}$ with $\sigma^*(V_{\fata})=V_{\fatb}$ 
  (according to \ref{cormot}\ref{cormotdiagrams}); 
  so the $\{v_b\}$ indeed form a basis of $V_{[\fata]}$.\\[2 mm]
  To see that $v_b$ does not depend on the choice of $\sigma_b$, let $\sigma\in S_n$ be another
  permutation with $\sigma^{-1}(\fata)=\fatb$.  
  It follows that $\fatb$ is fix under $\omega:=\sigma_b^{-1}\sigma$.
  Without loss of generality, we can assume that  
  $\omega=[1,\ldots,k_1]\cdot[k_1+1,\ldots,k_2]\cdot\ldots\cdot[k_{s-1}+1,\ldots,k_s]$.
  Because $\omega$ fixes $\fatb$, we have
  $b_1=b_2=\ldots=b_{k_1},\ b_{k_1+1}=\ldots=b_{k_2},\ b_{k_{s-1}+1}=\ldots=b_{k_s}$.
  Now express each cycle in this product as a product of transpositions that only permutes
  elements which are also permuted by the respective cycle. 
  Then these transpositions fix $\fatb$ as well, i.e. we can express $\omega$ as a product of
  transpositions that fix $\fatb$.
  Now \ref{lemmaminuseins} implies that $\omega^*\vert V_{\fatb}$ is multiplication by 
  $\sign\,(\omega)$, and the independance follows:
  \begin{multline}\label{eqindependant}
    \sign\,(\sigma)\cdot\sigma^*v
    =\sign\,(\sigma_b\omega)\cdot(\sigma_b\omega)^*v
    =\sign\,(\sigma_b\omega)\cdot\omega^*\sigma_b^*v\\
    =\sign\,(\sigma_b\omega)\cdot\sign\,(\omega)\cdot\sigma_b^*v
    =\sign\,(\sigma_b\omega^2)\cdot\sigma_b^*v
    =\sign\,(\sigma_b)\cdot\sigma_b^*v
    =v_b.
  \end{multline}
  To prove \eqref{eqvb}, take any $\sigma\in S_n$ and $\fatb\in[\fata]$. 
  We have
  $v_{\sigma^{-1}b}=v_{\sigma^{-1}\sigma_b^{-1}a}=
  v_{(\sigma_b\sigma)^{-1}a}\stackrel{\eqref{eqindependant}}{=}
  \sign\,(\sigma_b\sigma)\cdot(\sigma_b\sigma)^*v$, so that we get:
  \begin{multline*}
    \sigma^*v_b=\sigma^*(\sign\,(\sigma_b)\cdot\sigma_b^*v)
    =\sign\,(\sigma_b)\cdot\sigma^*\sigma_b^*v\\
    =\sign\,(\sigma_b)\cdot\underbrace{\sign\,(\sigma)\cdot\sign\,(\sigma)}_{=1}\cdot
      (\sigma_b\sigma)^*v
    =\sign\,(\sigma)\cdot\underbrace{\sign\,(\sigma_b\sigma)\cdot
      (\sigma_b\sigma)^*v}_{=v_{(\sigma_b\sigma)^{-1}a}=v_{\sigma^{-1}b}}
    =\sign\,(\sigma)\cdot v_{\sigma^{-1}b}.
  \end{multline*}
  From this, \eqref{eqvchi} and \eqref{eqdual}, equation \eqref{eqvb} follows immediately.
\end{proof}

\vspace{\abstand}

\noindent
Putting everything together, we get our final result:

\begin{thm}\label{thmzeta}
  Let $m\in\N^{\geq 3}$ and $n\in\N^{\geq 2}$ be natural numbers,
  let $k=\F_q$ be a finite field with $\text{char}\,{k}>m$ and $\mu_m\subseteq k^\times$,
  and let $Q$ be a $\FPZ/\FPz$-form of the Fermat equation $\Pmn$ 
  with associated hypersurface $Y$.\\[3 mm]
  Then because of \ref{thmclassification}\ref{efpz}, $Q\cong\Pmn\{L,x\}$ for
  $L=\prod_{i=1}^r\F_{q^{n_i}}$ with $\sum n_i=n$ 
  and $x=(x_1,\ldots,x_r)$ with $x_i\in{\F_{q^{n_i}}^\times}$. 
  According to \ref{satzfrobenius2},
  the associated cohomology class $(L,x)$ in $\Hc(G_k,\kranz)$ is represented by a
  1-cocycle $b$ that on the arithmetic Frobenius $f$ is given by
  $b_f=\prod_{i=1}^r\bigl[\bigl(y_i^{q^{n_i}-1},1,\ldots,1\bigr)\cdot z_{n_i}\bigr]\in\kranz$
  if $y_i$ denotes an $m$-th root of $x_i$ in $\bar{\F}_q$.
  Write $b_f$ in the form $(\zeta_1,\ldots,\zeta_n)\cdot\sigma$ with $\zeta_i\in\mum$ 
  and $\sigma\in S_n$.\\[4 mm]
  Then we know from \eqref{eqopfin}, 
  that the action of the gemetric Frobenius $F_Y$ on
  $V$ is given as
  \[
    F_Y^*
    =\Het\left(\op{\bar{f}^{-1}}{b_f}\right)\circ\FXs
    \stackrel{\mum\subseteq\F_q^\times}{=}\Het(b_f)\circ\FXs
    =\bigl((\zeta_1,\ldots,\zeta_n)\cdot\sigma\bigr)^*\circ\FXs,
  \]
  and it follows from \ref{satzjacobi}, \ref{satzkohom} and \ref{corbasis} that:
  \[\fbox{$\displaystyle
    F_Y^*\vert V_{\text{prim}}=
    \bigoplus_{[\fata]\in S_n\backslash A^m_n}
    \left\{
      V_{[\fata]}\xrightarrow{\cdot{\mathcal J}(\fata)}
      V_{[\fata]} 
      \xrightarrow{v_b\mapsto
        (\zeta_1^{b_{\sigma(1)}}\cdot\ldots\cdot\zeta_n^{b_{\sigma(n)}})\cdot\sign\,(\sigma)
        \cdot v_{\sigma^{-1}b}}
      V_{[\fata]}
    \right\}
  $},\]
  and finally from \ref{satzzeta} and \ref{satzkohom} we get
  \[\fbox{$\displaystyle
    \zeta(Q,t)=
    \text{det}\,\left(1-F_Y^*t\,\vert\,V_{\text{prim}}\right)
    ^{\left[(-1)^{n+1}\right]}\prod_{i\in\{0,\ldots,n-2\}}
    \frac{1}{1-q^it}
  $}.\]
\end{thm}

\vspace{\abstand}

\begin{bsp}
  Let $m=3$, $n=6$, $k=\F_7$ and $b\in\Hc(G_{\F_7},\mu_3\int S_6)$
  be as in \ref{bspintro} and \ref{bspjacobi}.
  We want to compute the zeta function of $Q:=P^3_6\{b\}$.
  According to \ref{bspintro}, on the arithmetic Frobenius $f$
  the 1-cocycle $b$ has value $(2,1,1,1,4,1)\cdot[1234][56]$:
  \[
    \begin{array}{cl}
      & \text{det}\,\left(1-F_Y^*t\,\vert\,V_{\text{prim}}\right) \\[4 mm]
      = & \text{det}\,(1-f^*F^*t\,|\,V_{\{1,1,1,1,1,1\}})
        \cdot\text{det}\,(1-f^*F^*t\,|\,V_{[1,1,1,2,2,2]}) \\
      & \hfill\cdot\text{det}\,(1-f^*F^*t\,|\,V_{\{2,2,2,2,2,2\}}) \\[4 mm]
      \stackrel{\ref{bspjacobi}}{=} & (1-7^2t)\cdot[1-(-56-21\zeta)t]\cdot[1-(49t)^4]^4\\
      & \hfill\cdot[1-(49t)^2]^2[1-(-35+21\zeta)t] \\[4 mm]
      = & (1+91t+7^4t^2)(1-7^2t)(1-7^4t^2)^2(1-7^8t^4)^4.
    \end{array}
  \]
  So \ref{thmzeta} gives
  \begin{multline*}
    \zeta(P,t)=
    \frac{1}{1-t}\cdot\frac{1}{1-7t}\cdot\frac{1}{(1+91t+7^4t^2)(1-7^2t)(1-7^4t^2)^2(1-7^8t^4)^4} \\
    \cdot\frac{1}{1-7^3t}\cdot\frac{1}{1-7^4t},
  \end{multline*}
  and by taking the logarithm, we finally get
  \[
    \zeta(P,t)=\text{exp}\,\left(
    \frac{2710}{1}t+\frac{5897984}{2}t^2+\frac{13881660703}{3}t^3+\frac{33246893493864}{4}t^4+\ldots
    \right).
  \]
\end{bsp}

\vspace{\abstand}

%%%%%%%%%%%%%%%%%%%%%%%%%%%%%%%%%%%%%%%%%%%%%%%%%%%%%%%%%%%%%%%%%%%%%%%%%%%%%%%%%%%%%%%%%%%%%%%%%%%%%%%%%%%%%%
%%%%%%%%%%%%%%%%%%%%%%%%%%%%%%%%%%%%%%%%%%%%%%%%%%%%%%%%%%%%%%%%%%%%%%%%%%%%%%%%%%%%%%%%%%%%%%%%%%%%%%%%%%%%%%
%%%%%%%%%%%%%%%%%%%%%%%%%%%%%%%%%%%%%%%%%%%%%%%%%%%%%%%%%%%%%%%%%%%%%%%%%%%%%%%%%%%%%%%%%%%%%%%%%%%%%%%%%%%%%%
%%%%%%%%%%%%%%%%%%%%%%%%%%%%%%%%%%%%%%%%%%%%%%%%%%%%%%%%%%%%%%%%%%%%%%%%%%%%%%%%%%%%%%%%%%%%%%%%%%%%%%%%%%%%%%
%%%%%%%%%%%%%%%%%%%%%%%%%%%%%%%%%%%%%%%%%%%%%%%%%%%%%%%%%%%%%%%%%%%%%%%%%%%%%%%%%%%%%%%%%%%%%%%%%%%%%%%%%%%%%%
%%%%%%%%%%%%%%%%%%%%%%%%%%%%%%%%%%%%%%%%%%%%%%%%%%%%%%%%%%%%%%%%%%%%%%%%%%%%%%%%%%%%%%%%%%%%%%%%%%%%%%%%%%%%%%
%%%%%%%%%%%%%%%%%%%%%%%%%%%%%%%%%%%%%%%%%%%%%%%%%%%%%%%%%%%%%%%%%%%%%%%%%%%%%%%%%%%%%%%%%%%%%%%%%%%%%%%%%%%%%%
%%%%%%%%%%%%%%%%%%%%%%%%%%%%%%%%%%%%%%%%%%%%%%%%%%%%%%%%%%%%%%%%%%%%%%%%%%%%%%%%%%%%%%%%%%%%%%%%%%%%%%%%%%%%%%
%%%%%%%%%%%%%%%%%%%%%%%%%%%%%%%%%%%%%%%%%%%%%%%%%%%%%%%%%%%%%%%%%%%%%%%%%%%%%%%%%%%%%%%%%%%%%%%%%%%%%%%%%%%%%%

\section{Remarks and open problems}

\vspace{\abstand}

Let us keep the notation of the last chapter but drop the assumption that
$k$ contains the $m$-th roots of unity, thus looking at a case 
which has not been studied in any of the above mentioned
papers by Weil, Deligne, Shioda or Gouv\^{e}a/ Yui.

In this situation, 
it is no longer true that the one-dimensional $\Q_l$-vector-spaces $V_{\fata}$ are invariant under
the Frobenius action. Indeed, in \cite{larsphd} we prove the following proposition:

\vspace{\abstand}

\begin{satz}\label{satzfrobenius}
  We have an action of $\zmzs$ on $A^m_n$, which is defined by
  $t\cdot(a_1,\ldots,a_n):=(ta_1,\ldots,ta_n)$. Because of $m>p$, we have $q\in\zmzs$,
  and for $\fata\in A^m_n$ we therefore can consider 
  $q\fata\in A^m_n$,
  $\langle\fata\rangle:=\{q^n\fata\vert n\geq 0\}$ 
  and $V_{\langle a\rangle}:=\bigoplus_{b\in\langle a\rangle}V_b$.\\[3 mm]
  The geometric Frobenius $F^*$ maps $V_a$ to $V_{q a}$.
  In particular, $V_{\langle a\rangle}$ is invariant under $F^*$ and therefore
  a $\Q_l$-$G_k$-representation, so that we get a decomposition of
  $V_{\text{prim}}$ in $\Qlk$ as
  \[
    V_{\text{prim}}=
    \bigoplus_{\textstyle\langle\fata\rangle\in\zmzs\backslash A^m_n}V_{\textstyle\langle\fata\rangle}.
  \]
  Put $d:=\text{gcd}\,(m,a_1,\ldots,a_n)$, $m':=m/d$,
  $\fatap:=(a_1/d,\ldots,a_n/d)\in A^{m'}_n$,
  and let $e$ denote the order of $q$ in $\zmzs$
  (so that $\F_{q^e}$ is the smallest field of characteristic $p$ that contains the 
  $m'$-th roots of unity).\\[3 mm]
  Then $\langle\fata\rangle=\{q^i\fata\,\vert\,i\in\{0,\ldots,e-1\}\}$
  with all the $q^i\fata$ being distinct.
  In particular, we have $e=\#\langle\fata\rangle$.\\[4 mm]
  Choose any $v\in V_{\fata}\setminus\{0\}$,
  and put $v_i:=(F^*)^iv\in V_{q^i a}\setminus\{0\}$ for $i\in{\{0,\ldots,{e-1}\}}$. 
  Then the matrix of $F^*\vert V_{\langle a\rangle}$ 
  with respect to the basis $\{v=v_0,v_1,\ldots,v_{e-1}\}$ is:
  \[
    \left(\begin{tabular}{ccc|c}
      $0$ & $\cdots$ & $0$ & ${\mathcal J}_{\F_{q^e}}^{m'}(\fatap)$ \\
      \hline
      $1$ &          & $0$ & $0$      \\
          & $\ddots$ &     & $\vdots$ \\
      $0$ &          & $1$ & $0$
    \end{tabular}\right)
  \]
\end{satz}

\vspace{\abstand}

This proposition enables us to explicitly compute the Frobenius action for forms of the Fermat hypersurface,
provided there are no
$\fata\in A^m_n$, $\sigma\in S_n\setminus\{\text{id}\}$ and $t\in\zmzs\setminus\{1\}$ 
with $t\fata=\op{\sigma}{\fata}$,
because then we can choose the bases of $V_{\text{prim}}$ described in
\ref{satzfrobenius} and \ref{corbasis} independantly. 

Such a ``good'' case is for example given by $m=n=3$, 
because in this case $A^3_3=\{(1,1,1),(2,2,2)\}$, 
$S_3$ acts trivially, and $(\Z/3\Z)^\times$ permutes the two characters.
This means we can compute Frobenius action and zeta function for arbitrary forms of the 
``Fermat cubic curve'' $x^3+y^3+z^3=0$ over any finite field of characteristic $\geq 5$.

In contrast to that, ``bad'' cases are for example given by $m=3$, $n=2$ or $m=3$, $n=4$:
\[
  2\cdot(1,2)=\op{(12)}{(1,2)} 
  \text{\ \ \ \resp\ \ \ }
  2\cdot(1,1,2,2)=\op{(13)(24)}{(1,1,2,2)}.
\]
Also in these two cases, we actually succeeded to compute the Frobenius action
by a very explicit combinatorial computation
(as explained in \cite{larsphd}), but to solve the problem in general, 
for any ``bad'' triple $(\fata,\sigma,t)$ as above, one has to compute the eigenvalue of
the endomorphism $(F\circ\sigma)^*\vert V_{\fata}$,
and we don't know how to do that yet.\\

This problem may also be connected to the Hodge conjecture and the Tate conjecture for Fermat varieties:

As we already mentioned in the introduction, 
Shioda proved thouse in a lot of cases (\cite{shioda79b}, \cite{shioda79}) 
and also identified the ``smallest'' case in which his method does not work in \cite{shioda81}: 
For $m=25$, $n=6$ and $p\equiv 1\pmod{m}$, the cohomology classes in $\fata:=(1,6,11,16,21,20)\in A^{25}_6$
are invariant under the Frobenius action, i.e. they are so-called ``Tate-cycles''
(\resp ``Hodge-cycles'' if one instead considers the same Fermat variety and $\fata$ defined over $\C$),
but Shioda's method
fails to prove that they are algebraic which they should be according to the Tate (\resp Hodge) conjecture.

Now we can not help noticing that $\fata$ is part of a ``bad'' triple in the above sense:
\begin{multline*}
  6\cdot(1,6,11,16,21,20)=(6,36,66,96,126,120)\\
  =(6,11,16,21,1,20)=\op{(12345)}{(1,6,11,16,21,20)},
\end{multline*}
suggesting that there may be a connection between the conjectures in this case and the
action of $\bigl(F\circ(12345)\bigr)^*$ on $V_{\fata}$ for $q\equiv 6\pmod{25}$. 

\vspace{\abstand}

%%%%%%%%%%%%%%%%%%%%%%%%%%%%%%%%%%%%%%%%%%%%%%%%%%%%%%%%%%%%%%%%%%%%%%%%%%%%%%%%%%%%%%%%%%%%%%%%%%%%%%%%%%%%%%
%%%%%%%%%%%%%%%%%%%%%%%%%%%%%%%%%%%%%%%%%%%%%%%%%%%%%%%%%%%%%%%%%%%%%%%%%%%%%%%%%%%%%%%%%%%%%%%%%%%%%%%%%%%%%%
%%%%%%%%%%%%%%%%%%%%%%%%%%%%%%%%%%%%%%%%%%%%%%%%%%%%%%%%%%%%%%%%%%%%%%%%%%%%%%%%%%%%%%%%%%%%%%%%%%%%%%%%%%%%%%
%%%%%%%%%%%%%%%%%%%%%%%%%%%%%%%%%%%%%%%%%%%%%%%%%%%%%%%%%%%%%%%%%%%%%%%%%%%%%%%%%%%%%%%%%%%%%%%%%%%%%%%%%%%%%%
%%%%%%%%%%%%%%%%%%%%%%%%%%%%%%%%%%%%%%%%%%%%%%%%%%%%%%%%%%%%%%%%%%%%%%%%%%%%%%%%%%%%%%%%%%%%%%%%%%%%%%%%%%%%%%
%%%%%%%%%%%%%%%%%%%%%%%%%%%%%%%%%%%%%%%%%%%%%%%%%%%%%%%%%%%%%%%%%%%%%%%%%%%%%%%%%%%%%%%%%%%%%%%%%%%%%%%%%%%%%%
%%%%%%%%%%%%%%%%%%%%%%%%%%%%%%%%%%%%%%%%%%%%%%%%%%%%%%%%%%%%%%%%%%%%%%%%%%%%%%%%%%%%%%%%%%%%%%%%%%%%%%%%%%%%%%
%%%%%%%%%%%%%%%%%%%%%%%%%%%%%%%%%%%%%%%%%%%%%%%%%%%%%%%%%%%%%%%%%%%%%%%%%%%%%%%%%%%%%%%%%%%%%%%%%%%%%%%%%%%%%%
%%%%%%%%%%%%%%%%%%%%%%%%%%%%%%%%%%%%%%%%%%%%%%%%%%%%%%%%%%%%%%%%%%%%%%%%%%%%%%%%%%%%%%%%%%%%%%%%%%%%%%%%%%%%%%

\bibliographystyle{alpha}
\bibliography{Literatur}

\end{document}